\documentclass[preprint]{elsarticle}

\usepackage{graphicx,natbib}
\usepackage{amsmath,amsfonts,amssymb,booktabs}
\usepackage{graphicx,epsfig,color,url,lscape}
\usepackage{subcaption,lscape,longtable,threeparttablex}
\usepackage{tikz}
\usepackage{float}
\usepackage{multirow}

\def\Cf{\mathcal{C}_{\mathcal{F}}}
\def\Ct{\mathcal{C}_{\mathcal{T}}}

\journal{}
\date{}
\begin{document}
\date{}
\begin{frontmatter}
\title{Parallel drone scheduling vehicle routing problems with collective drones}

\author{Roberto Montemanni, Mauro Dell'Amico, Andrea Corsini}


\address{Department of Sciences and Methods for Engineering, University of Modena and Reggio Emilia, Via Amendola 2, 42122 Reggio Emilia, Italy}

\begin{abstract}
We study last-mile delivery problems where trucks and drones collaborate to deliver goods to final customers. In particular, we focus on problem settings where either a single truck or a fleet with several homogeneous trucks work in parallel to drones, and drones have the capability of collaborating for delivering missions. This cooperative behaviour of the drones, which are able to connect to each other and work together for some delivery tasks, enhance their potential, since connected drone has increased lifting capabilities and can fly at higher speed, overcoming the main limitations of the setting where the drones can only work independently.

In this work, we contribute a Constraint Programming model and a valid inequality for the version of the problem with one truck, namely the \emph{Parallel Drone Scheduling Traveling Salesman Problem with Collective Drones} and we introduce for the first time the variant with multiple trucks, called the \emph{Parallel Drone Scheduling Vehicle Routing Problem with Collective Drones}. For the latter variant, we propose two Constraint Programming models and a Mixed Integer Linear Programming model.

An extensive experimental campaign leads to state-of-the-art results for the problem with one truck and some understanding of the presented models' behaviour on the version with multiple trucks. Some insights about future research are finally discussed.
\end{abstract}
\begin{keyword}

Parallel Drone Scheduling Vehicle Routing Problems with Cooperative Drones \sep Constraint Programming \sep Mixed Integer Linear Programming \sep Parallel Drone Scheduling Traveling Salesman Problems with Cooperative Drones
\end{keyword}
\end{frontmatter}

\section{Introduction}
The employment of drones in last-mile delivery is considered extremely strategic for the near future by leading distribution operators. They face a continuously increasing volume of parcels to handle, mainly generated by e-commerce (Statista, \cite{sta}). Considering that drones are light-weighted and use low-emission electric motors, that they do not have to move along the road network but can fly approximately in straight lines, and that they are not affected by road traffic congestions, their adoption for deliveries could lead to advantages for the companies (operational costs reduction), for the customers (faster deliveries) and for the whole society (sustainability).
Forbes \cite{for} refers to the heavy interest in drone technology as the ``Drone Explosion". The authors of \cite{BCG} forecast that autonomous vehicles will deliver about 80\% of all parcels in the upcoming decade.
In this work, we analyze a transition scenario where drones are used in conjunction with trucks for last-mile delivery.

Murray and Chu \cite{murray2015flying} introduced the idea of a new routing problem in which a truck and a drone collaborate to make deliveries. The authors present two new prototypical models expanding from the traditional Traveling Salesman Problem (TSP) called the Flying Sidekick TSP (FSTSP) and the Parallel Drone Scheduling TSP (PDSTSP). In both cases, a truck and some drones collaborate to deliver parcels. In the former model, drones can be launched from the truck during its tour, while in the latter one, drones are only operated from the central depot, and the truck executes a traditional delivery tour. In the remainder of the paper, we will focus on the latter problem, addressing the interested reader, e.g., to \cite{amicobb} and \cite{novel} for details and solution strategies for the FSTSP.

More formally, in the PDSTSP there is a truck that can leave the depot, serve a set of customers, and goes back to the depot. In parallel, there is also a set of drones, and each one of them can leave the depot, serve a customer, and return to the depot before serving other customers. Some of the customers cannot be served by the drones, either due to their location or the characteristic of their parcel. The objective of the optimization is to minimize the completion time of the last vehicle returning to the depot (or a cost function related to this) while serving all the customers.

A first Mixed Integer Linear Programming (MILP) model for the PDSTSP is proposed in \cite{murray2015flying} together with some simple heuristic methods. Another MILP model and the first metaheuristic method, based on a two steps strategy embedding a dynamic programming component, are discussed in \cite{mbiadou2018iterative}. Another two steps approach is presented in \cite{DMN} while a hybrid ant colony optimization metaheuristic is discussed in \cite{dinh2022} and a variable neighbor search one in \cite{lei2022}. In \cite{md23}, an effective constraint programming approach is proposed, which optimally solved all the benchmark instances previously adopted in the literature for both exact and heuristic methods.
Recently, in \cite{HA} another exact approach based on branch-and-cut was proposed, together with new benchmark problems.

Several PDSTSP variants are also introduced and studied in the literature, see e.g., \cite{ottooptimization} and \cite{pasha2022} for extensive surveys. We review herein only those extensions of the original problem that we find more relevant to the present study.

The recent work \cite{mbiadou2022} discusses the \emph{Parallel Drone Scheduling Multiple Traveling Salesman Problem}, which is a straightforward extension of the PDSTSP where multiple trucks are employed and the target is to minimize the time required to complete all the customer delivery. 
The authors propose a hybrid metaheuristic algorithm, a mixed integer linear model, and a branch-and-cut approach. The same problem is independently introduced also in \cite{raj2021}, where the authors propose three mixed integer linear programming models, together with a branch-and-price approach. A heuristic version of the branch-and-cut method is also introduced, aiming at solving the larger instances. A more realistic variation of the PDSTSP is introduced in \cite{nguyen2022}. In this version of the problem concepts such as capacity, load balancing, and decoupling of costs and times are taken into account. The authors propose a mixed integer linear programming model and a \emph{Ruin\&Recreate} metaheuristic for the problem. Constraint Programming methods for these variants of the PDSTSP employing several trucks, are discussed in \cite{mdvrp}, where convincing experimental results are also presented.

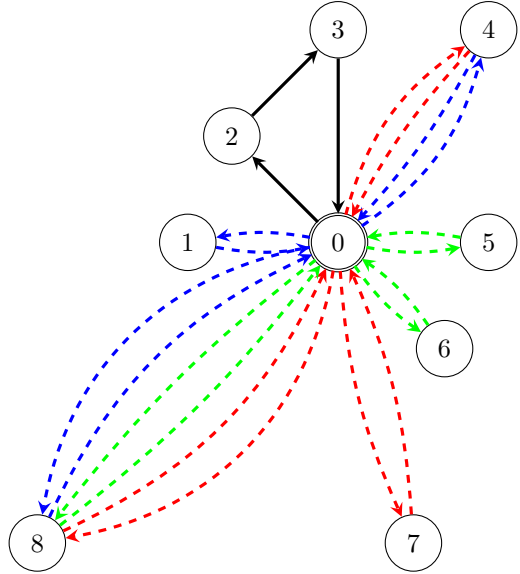
\begin{figure}[h!]
{
\begin{center}
\begin{tikzpicture}[node distance={2cm}, main/.style = {draw, circle}]
			\node[main,minimum size=0.75cm,double] (0) {0};
			\node[main,minimum size=0.75cm] (1) [left of=0] {1};
			\node[main,minimum size=0.75cm] (2) [above left of=0] {2};
			\node[main,minimum size=0.75cm] (3) [above right of=2] {3};
			\node[main,minimum size=0.75cm, label={[xshift=0.3cm,]}] (4) [right of=3] {4};
			\node[main,minimum size=0.75cm, label={[xshift=0.3cm,]}] (6) [below right of=0] {6};
			\node[main,minimum size=0.75cm, label={[xshift=0.3cm,]}] (7) [yshift=-4cm,xshift=1.cm] {7};	
			\node[main,minimum size=0.75cm, label={[xshift=-0.3cm,]}] (8) [yshift=-4cm,xshift=-4cm] {8};
			\node[main,minimum size=0.75cm, label={[xshift=0.3cm,]}] (5) [right of=0 ] {5};			
			\path [draw,->,>=stealth, line width=1.2] (0) ->  (2);
			\path [draw,->,>=stealth, line width=1.2] (2) -> (3);
			\path [draw,->,>=stealth, line width=1.2] (3) -> (0);
			\path [->,>=stealth, line width=1.2, color=blue,dashed] (0) edge [bend right=35] node[auto] {} (8);
			\path [->,>=stealth, line width=1.2, color=blue,dashed] (8) edge [bend right=-20] node[auto] {} (0);
			\path[->,>=stealth, line width=1.2, color=red,dashed] (0) edge [bend right=-35] node[auto] {} (8);
			\path[->,>=stealth, line width=1.2, color=red,dashed] (8) edge [bend right=20] node[auto] {} (0);
			\path[->,>=stealth, line width=1.2, color=green,dashed] (0) edge [bend right=7] node {} (8);
			\path[->,>=stealth, line width=1.2, color=green,dashed] (8) edge [bend right=7] node {} (0);
			\path[->,>=stealth, line width=1.2, color=green,dashed] (0) edge [bend right=10] node {} (5);
			\path[->,>=stealth, line width=1.2, color=green,dashed] (5) edge [bend right=10] node {} (0);
			\path[->,>=stealth, line width=1.2, color=green,dashed] (0) edge [bend right=10] node {} (6);
			\path[->,>=stealth, line width=1.2, color=green,dashed] (6) edge [bend right=10] node {} (0);
			\path[->,>=stealth, line width=1.2, color=blue,dashed] (0) edge [bend right=20] node {} (4);
			\path[->,>=stealth, line width=1.2, color=blue,dashed] (4) edge [bend right=-7] node {} (0);
			\path[->,>=stealth, line width=1.2, color=red,dashed] (0) edge [bend right=-20] node {} (4);
			\path[->,>=stealth, line width=1.2, color=red,dashed] (4) edge [bend right=7] node {} (0);
			\path[->,>=stealth, line width=1.2, color=blue,dashed] (0) edge [bend right=10] node {} (1);
			\path[->,>=stealth, line width=1.2, color=blue,dashed] (1) edge [bend right=10] node {} (0);
			\path[->,>=stealth, line width=1.2, color=red,dashed] (0) edge [bend right=10] node {} (7);
			\path[->,>=stealth, line width=1.2, color=red,dashed] (7) edge [bend right=10] node {} (0);
		\end{tikzpicture}
	\caption{Example of a PDSTSP-c instance. Node 0 is the depot, the other nodes are customers. Travel times are omitted for the sake of simplicity. The black {continuous} arcs represent the tour of the truck (0,\,2,\,3,\,0). The dashed  arcs depict the missions of the drones, each colour representing a different one. Note that some of the missions are carried out by multiple drones.}
	\label{figu}
\end{center}
}
\end{figure}

One common assumption in the literature on combined truck-drone delivery models has been a linear battery consumption for drones, leading to fixed operation ranges and carrying capacities. Recently, power consumption models with more realistic settings have been presented, {e.g., in \cite{raj2021}, \cite{raj20}, and \cite{liu17}, where the impact of a drone's power consumption is analyzed as a function of both speed and payload.} A review of drone energy consumption models is also available in \cite{zha21}.  In \cite{pac16}, a novel method based on advanced power consumption models and using a so-called ``Collective Drone'' (c-drone) is introduced. In these settings, multiple drones may be coupled together to aerially transport {items of large size and weight}. By sharing resources, such as power and operating instructions, a collective drone might outperform a single drone to operate more efficiently. The authors of \cite{viet} joined these ideas to come out with an innovative problem, called the PDSTSP-c, where c stands for \emph{collective}. In this problem, a realistic model is used to calculate the endurance and capacity of groups of drones working together to carry out tasks. The authors are able to pre-compute the optimal speed to carry out a certain delivery with different (smaller or larger) formations of drones. Based on these calculations, they propose a mixed integer programming model and a \emph{Ruin\&Recreate} metaheuristic for the newly introduced problem.
An example of a PDSTSP-c instance is provided in Figure \ref{figu}.

The contributions of the present paper are as follows:
\begin{itemize}
\item A new Constraints Programming model for the PDSTSP-c is {introduced}, together with a valid inequality. Experimental results show the great potential of the new model;
\item The PDSVRP-c problem is firstly introduced, where the settings of the PDSTSP-c are kept, but a fleet of vehicles is available instead of a single truck;
\item Two new Constraints Programming models and a Mixed Integer Programming model for the PDSVRP-c are {introduced} and validated through some experimental tests.
\end{itemize}

The remainder of the paper is organized as follows. In Section \ref{tsp}, the PDSTSP-c is formally described and a new Constraint Programming model is introduced, together with a new valid inequality. Section \ref{vrp} firstly introduces the PDSVRP-c as an extension of the former problem. Two Constraint Programming models and a Mixed Integer Linear Programming model are presented. Section \ref{resu} presents experimental results for the two problems, while conclusions are drawn in Section \ref{conc}.

\section{The Parallel Drone Traveling Salesman Problem with Collective Drones}\label{tsp}
In this section, we formally describe the PDSTSP-c, as originally introduced in \cite{viet}, and we present a new Constraint Programming model. We start from the single-vehicle problem because its description helps in introducing the multiple-vehicle version.

\subsection{Problem Description}\label{tspp}
Given a complete graph $G(V, E)$ with the set of vertices $V = \{0, 1, \dots , n\}$, with vertex 0 representing the depot, and the remaining vertices being associated with the customers (set $C = V \setminus \{ 0 \}$). Each customer $i$ requests delivery of a parcel of weight $w_i$ from the depot. The fleet of vehicles available for deliveries is a driver-operated delivery truck, with unlimited range and capacity, and a set $D$ of $m$ homogeneous drones that are based at the depot and equipped with batteries of given capacity (a fresh battery is installed before each mission). The truck performs its task within a single tour, beginning from the depot, traversing through all assigned {customers}, and returning to the depot. The truck travel times between pairs of vertices $i, j \in V$ is given as $t_{ij}$. {Matrix $[t_{ij}]$ satisfies the triangular inequality $t_{ik}\leq t_{ij}+t_{jk}, i,j,k\in V$.} The drones have to perform back-and-forth trips between the depot and the customers' locations to deliver the parcels. Travel times and ranges of drone missions depend on factors such as the number of drones cooperating and the traveling speed. Given a customer $i$ and a number $k$ of drones executing the mission, it is possible to pre-calculate the optimal speed and consequently the total travel time $\tau^k_i$ {for the back-and-forth trip}. When it is not possible to service a customer $i$ for some values of $k$, then $\tau^k_i$ is set to $+\infty$. We group such customers in set $\Ct \subsetneq \mathcal{C}$.
{Instead, let $\Cf = \mathcal{C} \setminus \Ct$ denote the (sub)set of customers that may be served with some drones' configuration, and let $q_j$ and $p_j$ be the minimum and maximum number of drones to serve a customer $j \in \Cf$ (observe that in the cited model the number of drones that can serve a customer  always define an interval). }
We adopted the realistic model described in \cite{viet} for the calculation of their travel times, and we refer the interested reader to this paper for full details.

Note that a main difficulty of the problem is that once $k$ drones collaborate for a delivery mission, strict synchronization constraints must be fulfilled. The objective of the PDSTSP-c is to find {a truck tour, drone-customer assignments, and drones scheduling} that minimize the makespan (i.e., the maximum completion time at which all vehicles are back at the depot after completing their services) while fulfilling all the constraints and conditions listed above.

\subsection{A Constraint Programming model}\label{tspmodel}
The Constraint Programming model we present is based on the Google-OR CP-SAT solver \cite{ortools} and follows the ideas behind the Mixed Integer Linear Program described in \cite{viet}. In particular, drone missions are modeled through a flow. Changes have however been introduced to take full advantage of the characteristic of the solver used.

The variables used in the model are as follows:
\begin{itemize}
    \item $x_{ij}$: binary variables equal to 1 (\emph{true}) if edge $(i, j)$, {with $i, j \in \mathcal{C}$}, is traveled by the truck, 0 (\emph{false}) otherwise. {Whereas a loop $x_{jj}=1$ means that customer $j$ is served by drones, while $x_{jj}=0$ if it is served by the truck.}
    \item $z^k_{j}$: binary variables equal to 1 if  {customer $j \in \Cf$} is served by $k$ drones, 0 otherwise.
    \item $y_{ij}$: binary variables equal to 1 if customer $i$ is served right before customer $j$ within the schedule of any drone,  0 otherwise.
    \item $f_{ij} \in R^{+}$: continuous flow variables indicating number of drones serving customer $i$ right before customer $j$ in their schedule.
    \item $T_j \in R^{+}$: continuous variables representing the time at which the service of customer $j \in \mathcal{C}_{\mathcal{F}}$ is completed by the drones, $T_0$ is the availability  time for all vehicles at the depot.
    \item $\alpha \in R^{+}$: continuous variable denoting the completion time, by which all the carriers are back to the depot.
\end{itemize}
\begin{align}
  (CP1)&:  \ \ \ \min   \alpha & \label{1a}\\
 s.t. \ \ \        & \alpha \geq \sum_{i \in V} \sum_{j \in V, i \ne j} t_{ij} x_{ij} &  \label{2a}\\
       & \alpha \geq T_j  & j \in \Cf   \label{3a}\\
       & x_{jj} = \sum_{q_j \leq k \leq p_j} z^k_{j}  & j \in \Cf  \label{4a}\\
       &	\text{Circuit($x_{ij},$ with $ i,j \in V, j \neq i $ if $j \in \Ct$)} \!\!\!\!\!\!\!\!\!\!\!\!\!\!\!\!\!\!& \label{6a}\\
      & \sum_{ j \in \Cf } f_{0j} \leq m & \label{7a}\\
      & \sum_{ i \in \Cf \cup \{0\}, i \neq j } f_{ij} =   \sum_{q_j \leq k \leq p_j} k z_j^k&j \in \Cf \label{8a}\\
      & \sum_{ i \in \Cf \cup \{0\}, i \neq j } f_{ij} =   \sum_{ l \in \Cf \cup \{0\}, l \neq j } f_{jl} &j \in \Cf \cup \{0\} \label{9a}\\
      & f_{ij} \leq m y_{ij} & i, j \in \Cf \cup \{0\}, i \neq j  \label{10a}\\
      & y_{ij} \implies T_j \geq T_i + \sum_{ q_j \leq k \leq p_j} \tau^k_j z^k_j & i \in \Cf \cup \{0\}, j \in \Cf,i \neq j  \label{12a}\\
      & 0 \leq f_{ij} \leq m & i,j \in \Cf \cup \{0\}, i \neq j \label{13a} \\
	& x_{ij} \in \{0; 1\} & i, j \in V \label{15a}\\
	& z^k_{j} \in \{0; 1\} & j \in \Cf, q_j \leq k \leq p_j  \label{16a}\\
      & y_{ij} \in \{0; 1\} & i,j \in \Cf \cup \{0\}, i \neq j \label{17a} \\
      & T_{j} \geq 0 & j \in \Cf \cup \{0\} \label{18a}
\end{align}
Following the trivial objective function (\ref{1a}), the constraints have the following meaning.
Constraint (\ref{2a}) says that the total time $\alpha$ has to be greater than or equal to the time required by the truck tour.
Analogously, constraints (\ref{3a}) impose that $\alpha$ has to be greater than or equal to the completion time of the eventual drone mission to serve customer $j$.
Given the logic of the variables, constraints (\ref{4a}) state that each drone-eligible customer has to be visited either by the truck or by a group of drones;
Constraint (\ref{6a}) uses the CP-SAT {method} \emph{Circuit} \cite{ortools} to force a feasible truck tour, eventually skipping each customer $j$ for which $x_{jj}=1$. {Note that \emph{Circuit} is invoked giving it all the variables $x_{ij}$ with $i\neq j$ and all the variables $x_{jj}$ with $j\not\in \Ct$. This implies that $x_{jj}$ remains zero for all $j\in \Ct$, so all the customers that are not drone-eligible must be served by the truck. }
The Constraints (\ref{7a})-(\ref{9a}) model the operations and synchronization of the drones as a flow problem (see \cite{viet} for more detailed explanations):
Constraint (\ref{7a}) states the flow going out from node 0 {has to be less than or equal to $m$ (remind that each drone is represented as a unit of flow)};
Constraints (\ref{8a}) impose that if a customer $j$ is serviced by $k$ drones, than the flow entering node $j$ has to equal $k$;
Constraints (\ref{9a}) are classic conservation equalities, imposing that the flows entering and exiting a node must be equal.
{Constraints (\ref{10a}) activate the variables $y$ corresponding to arcs used by flows (variables $f$) which are necessary to calculate the completion time of drones.
Constraints (\ref{12a}) are active only if the variable $y_{ij}=1$ and {state that the synchronization constraint on arc (i, j) must be respected}.
This is achieved through the CP-SAT command \emph{OnlyEnforceIf} \cite{ortools}, which is indicated with $\implies$ in the model.
}
The remaining constraints (\ref{13a})-(\ref{18a}) define the domain of the variables.

\subsection{Valid inequality}\label{sva}
{The above basic model can be improved by means of the following new valid inequality}:
\begin{equation}
m \alpha \ge \sum_{j \in \Cf} \sum_{q_j \leq k \leq p_j} k \tau^k_j z^k_j \label{va}
\end{equation}
Inequality (\ref{va}) constraints $\alpha$ to be at least the total time spent on a mission by all the drones, divided by the number of drones used. Assuming that at optimality the drones will more or less balance the workload among them, makes the inequality effective. Inequality (\ref{va}) exposes directly mission times to $\alpha$, without filtering them through $T$ and $z$ variables, making the linear relaxation of the model much tighter.
Note however that (\ref{va}) fails to capture the time waited by the drones to synchronize with the others in case of multi-drone missions.
Finally, inequality (\ref{va}) will be valid also for all the models discussed in Section \ref{vrp} for the PDSVRP-c.

\section{The Parallel Drone Vehicle Routing Problem with Collective\\ Drones}\label{vrp}
In this section, we build upon Section \ref{tsp} and introduce the PDSVRP-c, a natural extension of the PDSTSP-c where multiple vehicles operate in parallel to the drones.
The problem is introduced in Section \ref{vrpp} while two models based on Constraint Programming are discussed in Sections \ref{vrp2} and \ref{vrp3}. {The first one is a 2-indices formulation based on the $CP1$  model of the previous section while the second one is a 3-indices formulation.}
{Section \ref{milp_vrp3} outlines a  MILP model to serve as a baseline.}

\subsection{Problem Description}\label{vrpp}
A formal definition of the PDSVRP-c can be proposed as a straightforward extension of the PDSTSP-c provided in Section \ref{tspp}. The difference is that now we have a fleet $S$ of $s$ trucks, with the same characteristics of the single truck employed for the PDSTSP-c: unlimited capacity, unlimited range, and same traveling speed. No concept of collaboration exists for the trucks and each customer has to be served either by one of the trucks or by drones.

Having a fleet does not change substantially the problem, but has an impact on the optimization since we now have to plan multiple tours and account for the mission time of each truck while calculating the completion time $\alpha$. We will see in the next sections two alternative {Constraint Programming} models and a Mixed Integer Linear Programming formulation.

\subsection{A 2-indices Constraint Programming model}\label{vrp2}
This model is the direct extension of that discussed in Section \ref{tspmodel} for the $CP1$ and delegates the Constraint Programming solver to handle the multiple truck tours. The variables remain the same, although now the $x$ can take the shape of multiple tours instead of a single one. {Another important difference is the definition of the variables $T_j$. In the $CP1$ model of Section \ref{tspmodel}, they are only related to the drones and represent the time in which the service of a customer is completed. Here, they are extended to the customers served by the trucks and represent the \emph{starting} time of the service of the truck to the customer. Formally we use the new variables $\overline T$ with the following meaning}
\begin{itemize}
\item {$\overline T_j \in R^{+}$: continuous variables representing the time at which the service of customer $j \in \mathcal{C}_{\mathcal{F}}$ is completed, if the customer is served by drones, or the service is started if the customer is served by a truck. $\overline T_0$ denotes the availabilty time at the depot for all vehicles. }
\end{itemize}

\begin{align}
  (CP2)&  \ \ \ \min   \alpha & \label{1b}\\
 s.t. \ \ \  & \alpha \geq \overline T_j  +  t_{j0} x_{j0} & j \in \mathcal{C} \label{2b}\\
    	& \alpha \geq \overline T_j  & j \in \Cf   \label{3b}\\
       & x_{jj} = \sum_{q_j \leq k \leq p_j} z^k_{j}  & j \in \Cf  \label{4b}\\
       &	\text{MultipleCircuit($x_{ij},$ with $ i,j \in V,  i \neq 0 \lor j \neq 0, j \neq i $ if $i \in {\Ct}$)} \!\!\!\!\!\!\!\!\!\!\!\!\!\!\!\!\!\!\!\!\!\!\!\!\!\!\!\!\!\!\!\!\!\!\!\!\!\!\!\!\!\!\!\!\!\!\!\!\!\!\!\!\!\!\!\!\!\!\!\!\!\!\!\!\!\!\!\!\!& \label{6b}\\
      & \sum_{ j \in \mathcal{C}} x_{0j} \leq s & \label{6ab}\\
	& x_{ij} \implies \overline T_j \geq \overline T_i + t_{ij}& i \in {V}, j \in \mathcal{C}, i \neq j  \label{6bb}\\	
      & \sum_{j \in \Cf} f_{0j} \leq m & \label{7b}\\
      & \sum_{ i \in \Cf \cup \{0\}, i \neq j } f_{ij} =   \sum_{q_j \leq k \leq p_j} k z_j^k&j \in \Cf \label{8b}\\
      & \sum_{ i \in \Cf \cup \{0\}, i \neq j } f_{ij} =   \sum_{ l \in \Cf \cup \{0\}, l \neq j } f_{jl} &j \in \Cf \cup \{0\} \label{9b}\\
      & f_{ij} \leq m y_{ij} & i, j \in \Cf \cup \{0\}, i \neq j  \label{10b}\\
      & y_{ij} \implies \overline T_j \geq \overline T_i + \sum_{ q_j \leq k \leq p_j} \tau^k_j z^k_j & \!\!\!\!\!\!\!\!\! i \in \Cf \cup \{0\}, j \in \Cf,i \neq j  \label{12b}\\
      & 0 \leq f_{ij} \leq m & i,j \in \Cf \cup \{0\}, i \neq j \label{13b} \\
	& x_{ij} \in \{0; 1\} & i, j \in {V} \label{15b}\\
	& z^k_{j} \in \{0; 1\} & j \in \Cf, q_j \leq k \leq p_j  \label{16ab}\\
      & y_{ij} \in \{0; 1\} & i,j \in \Cf \cup \{0\}, i \neq j \label{17b} \\
      & \overline T_{j} \geq 0 & j \in {V} \label{18b}
\end{align}

The constraints strictly follow the meaning already described for the $CP1$ model in Section \ref{tspmodel}. The only changes are as follows. Constraint (\ref{2b}) now defines $\alpha$ based on the time required by each truck to go back to the depot after visiting each of its assigned customers. {This constraint is valid since the travel times satisfy the triangular property, although it could be made valid also for the general case with the use of a \emph{OnlyEnforceIf} statement (see below).} Constraint (\ref{6b}) describes a set of circuits through the \emph{MultipleCircuit} command of CP-SAT (\cite{ortools})  to reflect we are now dealing with several tours instead of one. The new constraint (\ref{6ab}) forces the number of tours to be maximum $s$. Whereas the new constraints (\ref{6bb}) calculate the service start time for each customer visited by a truck (remember that $\implies$ indicates the \emph{OnlyEnforceIf}, which activates the constraint iif $x_{ij} = 1$).


\subsection{A 3-indices Constraint Programming model}\label{vrp3}

This model is another extension of the $CP1$ model in Section \ref{tspmodel} that uses $s$ separate sets of variables to describe the tours of the $s$ trucks.
All the variables remain the same, apart from the $x$s which are substituted by a set of variables $w$ such that $w^k_{ij} = 1$ if edge $(i, j)$ is traveled by truck $k \in S$, 0 otherwise. Note that $w^k_{jj}=1$ means that customer $j$ is not served by truck $k$, hence it is not part of its tour. In addition, $w_{00}^k=1$ means that truck $k$ is not operated in the solution. {Note that differently from model $CP2$, here all the loop variables for the truck are inserted when invoking the method \emph{Circuit} used to find a circuit for each truck, see \eqref{6c} below.}
Finally note that {the timing variables $T$ are the same used in model $CP1$ of Section \ref{tspmodel}}.

\begin{align}
  (CP3)&  \ \ \ \min   \alpha & \label{1c}\\
 s.t. \ \ \        & \alpha \geq \sum_{i \in V} \sum_{j \in V, i \ne j} t_{ij} w^k_{ij} &  k \in S \label{2c}\\
    	& \alpha \geq T_j  & j \in \Cf   \label{3c}\\
       & \sum_{k=1}^s (1-w^k_{jj}) = \sum_{q_j \leq k \leq p_j} z^k_{j}  & j \in \Cf  \label{4c}\\
       & \sum_{k=1}^s w^k_{jj} = s-1  & j \in \mathcal{C}_{\mathcal{T}}    \label{4ac}\\
       &	\text{Circuit($w^k_{ij},$ with $ i,j \in V$)} \!\!\!\!\!\!\!\!\!\!\!\!\!\!\!\!\!\! & k \in S \label{6c}\\
      &  w^k_{ij} \leq  1 - w^k_{00} & k \in S, i,j \in  \mathcal{C} \label{6ac}\\
      & \sum_{ j \in \Cf } f_{0j} \leq m & \label{7c}\\
      & \sum_{ i \in \Cf \cup \{0\}, i \neq j } f_{ij} =   \sum_{q_j \leq k \leq p_j} k z_j^k&j \in \Cf \label{8c}\\
      & \sum_{ i \in \Cf \cup \{0\}, i \neq j } f_{ij} =   \sum_{ l \in \Cf \cup \{0\}, l \neq j } f_{jl} &j \in \Cf \cup \{0\} \label{9c}\\
      & f_{ij} \leq m y_{ij} & i, j \in \Cf \cup \{0\}, i \neq j  \label{10c}\\
      & y_{ij} \implies T_j \geq T_i + \sum_{ q_j \leq k \leq p_j} \tau^k_j z^k_j & i \in \Cf \cup \{0\}, j \in \Cf,i \neq j  \label{12c}\\
      & 0 \leq f_{ij} \leq m & i,j \in \Cf \cup \{0\}, i \neq j \label{13c} \\
	& w^k_{ij} \in \{0; 1\} & k \in S, i, j \in V \label{15c}\\
	& z^k_{j} \in \{0; 1\} & j \in \Cf, q_j \leq k \leq p_j  \label{16ac}\\
      & y_{ij} \in \{0; 1\} & i,j \in \Cf \cup \{0\}, i \neq j \label{17c} \\
      & T_{j} \geq 0 & j \in \Cf \cup \{0\} \label{18c}
\end{align}

The constraints strictly follow the meaning already described for the $CP1$ model in Section \ref{tspmodel}. The changes reflect the presence of multiple trucks and are as follows. Inequalities (\ref{2c}) now constrain $\alpha$ to be equal to or larger than the length of the tour of each truck $k$. Equalities (\ref{4c}) now express that each drone-eligible customer has to be visited either by one of the trucks or the drones. The new constraints (\ref{4ac}) state that customers in $\Ct$ cannot be visited by drones, and have to be serviced by exactly one truck. Constraints (\ref{6c}) are now independently defined for each truck $k$, dropping the concept of giant-tour introduced for the $CP2$ model. The new technical constraint (\ref{6ac}) forces the circuit of a truck $k$ to be empty once the relative variable $w^k_{00}$ takes the value 1.


\subsection{A 3-indices Mixed Integer Linear Programming model}\label{milp_vrp3}

We finally present the Mixed Integer Linear {Programming (MILP)} formulation of the PDSVRP-c, which is based on the 3-indices $CP3$ model of Section \ref{vrp3}.
For sake of simplicity in the presentation of the model, we adopt  a new variable $u_j^k$ that takes value 1 if $k \in S$ serves customer $j$, 0 otherwise. Note that this variable can be defined as $u_j^k=1-w^k_{jj}$ in the logic of the $CP3$ model, but in the {MILP model the loop  variables $w^k_{jj}$ are not used}.

\begin{align}
({MILP}) & \ \ \min \alpha & \label{milp:1} \\\
 s.t. \ \ \  & \alpha \ge \sum_{i \in V} \sum_{j \in V, i \neq j} t_{ij}w^k_{ij} & k \in S \label{milp:2} \\
    & \alpha \ge T_j  & j \in \Cf \label{milp:3}\\
    & \sum_{k \in S} u^k_j + \sum_{q_j\le k\le p_j} z_j^k = 1 & j \in \Cf \label{milp:4} \\
    & \sum_{k \in S} u_j^k = 1 & j \in \mathcal{C}_{\mathcal{T}} \label{milp:5}\\
    & u_j^k \le u_0^k  &  j \in \mathcal{C}, k \in S  \label{milp:6}\\
    & \sum_{i \in {V}, i \neq j} w^k_{ij} + \sum_{l \in {V}, l \neq j} w^k_{jl} = 2u_j^h & j \in {V}, k \in S \label{milp:7} \\
    & \sum_{i, j \in H, i \neq j} w^k_{ij} \le |H| - 1 &  {H \subseteq\mathcal{C}}, k \in S \label{milp:SEC}\\
    & \frac{f_{ij}}{m} \leq y_{ij} \leq f_{ij} & i, j \in \Cf \cup \{0\}, i \neq j \label{milp:8} \\
    & T_j + M(1 - y_{ij}) \geq T_i + \sum_{ q_j \leq k \leq p_j} \tau^k_j z^k_j & \!\!\!\! i \in \Cf \cup \{0\}, j \in \Cf,i \neq j  \label{milp:9}\\
    & 0 \leq f_{ij} \leq m & i,j \in \mathcal{C}_{\mathcal{F}}, i \neq j \label{milp:10} \\
    & w^k_{ij} \in \{0; 1\} & k \in S, i, j \in {V}, i \neq j \label{milp:10a}\\
    & z^k_{j} \in \{0; 1\} & j \in \Cf, q_j \leq k \leq p_j \label{milp:11}\\
    & {u^k_{j} \in \{0; 1\}} & {j \in \mathcal{C},  k\in S} \label{milp:11bis}\\
    & y_{ij} \in \{0; 1\} & i,j \in \Cf \cup \{0\}, i \neq j \label{milp:12} \\
    & T_{j} \geq 0 & j \in \Cf \cup \{0\}   \label{milp:13}
\end{align}

The model minimizes the time to serve all the customers (\ref{milp:1}). Constraints (\ref{milp:2}) force $\alpha$ to be larger than any tour of the trucks, and inequalities (\ref{milp:3}) guarantee that $\alpha$ is larger than {the completion time of} any drone's mission time.
Equations (\ref{milp:4}) assign customers from $\Cf$ to either a drone or a truck, while constraints (\ref{milp:5}) force the customers that can be visited only by a truck ($\mathcal{C}_{\mathcal{T}}$) to receive such a visit. Inequalities (\ref{milp:6}) impose that a customer can be visited by a truck only if it is in use. Equalities (\ref{milp:7}) are flow conservation constraints for the truck tours.
Inequalities (\ref{milp:SEC}) are subtour elimination constraints \cite{dfj} and guarantee that truck tours are { circuits including the depot}. Note that these constraints are exponential in number, depending on any possible  subset $H$ of {$\mathcal{C}$}. In our implementation, they will be generated dynamically as described in Section \ref{sepa} below.
Constraints (\ref{milp:8}) refers to the flow of drones and guarantee that the number of drones going from customer $i$ to $j$ must be lower than $m$, only if there is a flow from $i$ to $j$.
Inequalities (\ref{milp:9}) regulate {completion of } the service time for the customers visited by the drones.
Finally, constraints (\ref{milp:10})-(\ref{milp:13}) define the domain of the variables.

\subsubsection{Separation {of the subtour elimination constraints}} \label{sepa}

{The subtour elimination constraints  \eqref{milp:SEC} are dynamically added to the MILP as \emph{Lazy constraints} (available in the most popular MILP solvers).
Specifically, the solver starts by disregarding the constraints declared \emph{lazy} and once a feasible integer solution is found it invokes a user defined separation procedure.
In our case, since the solution on hand is integer, the separation is a simple $O(|E|)$ exploration of the graphs $G^k=(V,E^k)$ with $E^k=\{(i,j)\in E: w_{ij}^k=1\}$ to look for subtours not involving the depot.}
\section{Experimental Results}\label{resu}
All the models presented in previous sections have been coded in Python 3.11.2.
The Constraint Programming models of Sections \ref{tsp} and \ref{vrp} have been solved via the CP-SAT solver of Google OR-Tools 9.6 \cite{ortools} while the Mixed Integer Linear model of Section \ref{vrp} has been solved with Gurobi 10.0 \cite{gurobi}.

{The outcome of the experimental campaign is discussed in the remainder of this section and is organized according to the different problems tackled.
{Tables \ref{tspts}-\ref{tab:vrp_large5} report, for each instance: i) the instance name; ii) the lower bound eventually produced and the best heuristic solution ([LB, UB]);  iii) the computing time to find the best heuristic solution (Sec$_{\text{bst}}$); iv) the eventual computing time to prove optimality (Sec$_{\text{tot}}$); v) a final summary column (Best bounds) containing the current state-of-the-art results of each instance for easing future research.}

In addition, we use a dash whenever a result is not retrieved {or the time limit is reached}, and we mark in \textit{italics} the results of our models not matching nor improving best-known bounds while in \textbf{bold} those producing new best bounds.} Hardware configurations, solvers used, experimental settings and  time limits are finally reported in the notes of the tables for each approach.

\subsection{Benchmark Instances}
To evaluate the performance of the proposed models for both the PDSTSP-c and the new PDSVRP-c, we consider the instances originally introduced in \cite{viet} for the PDSTSP-c, and available at \url{http://orlab.com.vn/home/download}. The number $n$ of customers varies from 15 to 200 (first number of the instance name) and the instances are divided into small ($n \le 30$) and large ($n > 30$). The number $m$ of drones available varies in the range $[3, 6]$ for the small instances and $[5, 10]$ for the large ones. {The traveling distances for trucks are computed using  Manhattan distances and a speed of 30 km/h, while drones follow the Euclidean distance and the optimal travel times (rounded up to the nearest integer) are pre-calculated for each collaborative cluster of $k$ drones. The interested reader can find all the details of the instances in  \cite{viet}.}

For generating PDSVRP-c instances {we used the same set of benchmarks and added the number $s$  of trucks chosen } in the range $[2, 3]$ for small instances and in $[2, 5]$ for large instances.

\subsection{PDSTSP-c}

In this section, we aim at comparing the  results obtained by solving the $CP1$ model described in Section \ref{tspmodel}, with and without the valid inequality (\ref{va}). The results are summarized in Table \ref{tspts} for the small instances and in Table  \ref{tsptl} for the large ones. We compare $CP1$ with the methods introduced in \cite{viet}, namely a Mixed Integer Linear Programming (MILP) model solved with IBM CPLEX 12.1 \cite{cplex} and two versions of a \emph{Ruin\&Recreate} metaheuristic: {RnR fast} and {RnR}. Note that the results of the MILP model in \cite{viet} are only available for small instances and those reported for the \emph{Ruin\&Recreate} methods are the best over 30 runs. To fully understand the impact of inequality (\ref{va}), we also considered the MILP model described in Section \ref{milp_vrp3} for the PDSVRP-c and run it with $s=1$ (one truck only) as well as the inequality (\ref{va}). This method is run on small instances only, since a previous study (\cite{viet}) showed MILP models are not suitable for the large instances.

We are not aware of other existing methods to deal with this problem.

\begin{landscape}
{
\begin{center}
\begin{table}[h]
\caption{Experimental results on the PDSTSP-c. Small instances. }\label{tspts} 
\resizebox{1.6\textheight}{!}{
\begin{threeparttable}
\begin{tabular}{l  cr  cr  cr  crr  crr  crr  c}
\hline
& \multicolumn{2}{c}{RnR fast  \cite{viet}\tnote{a}} & \multicolumn{2}{c}{RnR  \cite{viet}\tnote{a}} & \multicolumn{2}{c}{MILP \cite{viet}\tnote{b} } & \multicolumn{3}{c}{$MILP$+(\ref{va})\tnote{c} } & \multicolumn{3}{c}{$CP1$\tnote{d}} & \multicolumn{3}{c}{$CP1$+ (\ref{va})\tnote{d}} & Best \\
Instance & UB & Sec$_{\text{bst}}$ & UB & Sec$_{\text{bst}}$ & [LB, UB] & Sec$_{\text{tot}}$ & [LB, UB] & Sec$_{\text{tot}}$ & Sec$_{\text{bst}}$ & [LB, UB] & Sec$_{\text{tot}}$ & Sec$_{\text{bst}}$  & [LB, UB] & Sec$_{\text{tot}}$ & Sec$_{\text{bst}}$ & bounds \\
\cmidrule(lr){1-1}\cmidrule(lr){2-3}\cmidrule(lr){4-5}\cmidrule(lr){6-7}\cmidrule(lr){8-10}\cmidrule(lr){11-13}\cmidrule(lr){14-16}\cmidrule(lr){17-17}
15-r-e 	& 92 & 0.32 & 92 & 0.95 		& [92, 92] & 8.65 		&[92, 92]	&	1215.79	&	850.00& [92, 92] & 0.84 & 0.44 & [92, 92] & 1.14 & 1.10 & 92 \\
15-rc-c 	& 44 & 0.43 & 44 & 1.46 		& [31.74, 44] & - 		&[\textbf{44},	44]	&	1837.95	&	1790.09& [40, 44 ]& - & 28.79 & [\textbf{44}, 44] & 12.00 & 11.94 & 44 \\
16-c-c 	& 60 & 0.52 & 60 & 2.14 		& [60, 60] & 2.61 		&[\textbf{60},	60]	&	3.66	&	3.63& [\textbf{60}, 60] & 2.95 & 2.90 & [60, 60] & 2.10 & 2.05 & 60 \\
16-r-e 	& 112 & 0.55 & 112 & 1.52 	& [112, 112] & 63.86 		&[98.20,	112]	&	 -	&	250.10& [112, 112] & 2.00 & 1.94 & [112, 112] & 1.61 & 1.57 & 112 \\
18-c-c 	& 56 & 0.42 & 56 & 1.86 		& [38.72, 56] & - 		&[\textbf{56},	56]	&	1258.06	&	50.89& [{44}, 56] & - & 302.76 & [\textbf{56}, 56] & 42.22 & 42.13 & 56 \\
18-r-e 	& 96 & 0.59 & 96 & 1.79 		& [86.94, 96] & - 		&[87.86,	96]	&	 -	&	2338.15& [{92}, 96 ]& - & 55.83 & [\textbf{96}, 96] & 4.70 & 4.63 & 96 \\
18-rc-c 	& 58 & 0.54 & 57 & 2.40 		& [37.78, 57] & - 		&[56,	57]	&	 -	&	3060.21& [{38}, \textit{60}] & - & 399.43 & [\textbf{57}, 57] & 1803.93 & 6.57 & 57 \\
19-c-c 	& 44 & 0.48 & 44 & 1.81 		& [28.06, 44] & - 		&[40.85,	44]	&	 -	&	3302.58& [{32}, 44] & - & 64.50 & [\textbf{44}, 44] & 24.91 & 11.00 & 44 \\
20-c-c 	& 43 & 0.60 & 43 & 2.54 		& [30.99, 44] & - 		&[39.42,	43]	&	 -	&	2374.04& [\textbf{40}, 43] & - & 102.82 & [\textbf{40}, 43] & - & 147.26 & [40, 43] \\
20-r-c 	& 64 & 0.43 & 64 & 1.91 		& [55.35, 64] & - 		&[61.80,	64]	&	 -	&	3326.72& [{56}, 64] & - & 276.15 & [\textbf{64}, 64] & 73.39 & 73.29 & 64 \\
20-r-e 	& 82 & 0.62 & 80 & 2.00 		& [62.59, 88] & - 		&[72.80,	82]	&	 -	&	3237.98& [{72}, 80] & - & 606.13 & [\textbf{80}, 80] & 38.11 & 38.00 & 80 \\
20-rc-c 	& 96 & 0.41 & 96 & 2.37 		& [96, 96] & 0.9 		&[96,	96]	&	460.30	&	1.14& [96, 96] & 6.48 & 6.42 & [96, 96] & 6.23 & 6.17 & 96 \\
20-rc-e 	& 100 & 0.46 & 100 & 1.09 	& [100, 100] & 88.78 	&[90,	100]	&	 -	&	292.98& [100, 100] & 6.92 & 6.87 & [100, 100] & 4.70 & 4.65 & 100 \\
21-c-c 	& 62 & 0.48 & 62 & 2.25 		& [41.80, 64] & - 			&[44,	64]	&	 -	&	2281.04& [\textit{36}, \textit{64}] & - & 18.02 & [\textbf{60}, 62] & - & 58.44 & [60, 62] \\
21-r-e 	& 85 & 0.59 & 85 & 1.74 		& [59.52, 100] & -	 	&[75.25,	88]	&	 -	&	3572.98& [\textit{49}, \textit{88}] & - & 1512.11 & [\textbf{85}, 85] & 940.55 & 128.94 & 85 \\
23-c-e 	& 80 & 0.60& 80 & 2.49  		& [58.15, 80] & - 		&[58.15,	80]	&	 -	&	2698.31& [\textbf{80}, 80] & 0.84 & 0.78 & [\textbf{80}, 80] & 1.31 & 1.25 & 80 \\
23-r-c 	& 88 & 0.42 & 88 & 1.83 		& [88, 88] & 3293.16 	&[84,	88]	&	 -	&	3042.55& [88, 88] & 218.07 & 217.97 & [88, 88] & 8.97 & 8.90 & 88 \\
24-c-e 	& 84 & 0.79 & 84 & 2.50 		& [78.4, 84] & - 			&[78.40,	84]	&	 -	&	3567.08& [\textbf{84}, 84] & 14.07 & 13.98 & [\textbf{84}, 84] & 11.05 & 10.97 & 84 \\
24-r-e 	& 112 & 0.52 & 112 & 1.57 	& [91.05, 112] & - 		&[101,	112]	&	 -	&	1819.74& [{108}, 112] & - & 4.76 & [\textbf{112}, 112] & 955.94 & 3.60 & 112 \\
24-rc-c 	& 72 & 0.73 & 71 & 4.06 		& [69.58, 88] & - 		&[69.58,	72]	&	 -	&	14.13& [\textit{68}, 71] & - & 2245.17 & [\textbf{70}, \textbf{70}] & 190.03 & 189.86 & 70 \\
25-c-c 	& 56 & 0.60 & 56 & 2.92 		& [37.33, 56] & - 		&[37.83, 56]	&	-	& 	694.26 &[\textit{35}, 56] & - & 764.50 & [\textbf{56}, 56] & 44.85 & 38.62 & 56 \\
25-r-e 	& 106 & 0.96 & 104 & 3.14 	& [76.11, 120] & - 		&[95.73,	108]	&	 -	&	3533.88& [\textit{58}, \textit{108}] & - & 60.58 & [\textbf{104}, 104] & 288.19 & 288.04 & 104 \\
25-rc-e 	& 92 & 0.71 & 92 & 2.25 		& [66.99, 100] & - 		&[83.60,	96]	&	 -	&	2636.83& [\textit{60}, \textit{97}] & - & 90.45 & [\textbf{92}, 92] & 113.76 & 113.63 & 92 \\
26-r-c 	& 103 & 0.53 & 103 & 2.58 	& [95.26, 128] & - 		&[100.18,	103]	&	 -	&	3409.09& [\textit{84}, \textit{104}] & - & 2746.86 & [\textbf{101}, 103] & - & 107.06 & [101, 103] \\
27-c-c 	& 84 & 0.49 & 84 & 2.18 		& [83.23, 84] & - 		&[64.72,	84]	&	 -	&	2672.84& [{84}, 84] & 1.91 & 1.85 & [\textbf{84}, 84] & 1.80 & 1.75 & 84 \\
27-c-e 	& 68 & 0.72 & 68 & 6.27 		& [42.04, 68] & - 		&[42.04,	68]	&	 -	&	1429.86& [\textit{31}, 68] & - & 1.76 & [\textbf{68}, 68] & 33.23 & 1.36 & 68 \\
27-rc-c 	& 100 & 0.77 & 100 & 5.30 	& [100, 100] & 721.25 	&[85.34,	100]	&	 -	&	2915.17& [100, 100] & 394.35 & 394.19 & [100, 100] & 71.17 & 71.06 & 100 \\
27-rc-e 	& 84 & 0.79 & 84 & 2.70 		& [59.52, 100] & - 		&[64.29,	88]	&	 -	&	1202.70& [\textit{42}, \textit{88}] & - & 2066.04 & [\textbf{84}, 84] & 576.56 & 54.8 & 84 \\
29-rc-e 	& 116 & 0.72 & 116 & 1.69 	& [97.71, 124] & - 		&[109.75,	116]	&	 -	&	2214.71& [\textbf{116}, 116] & 654.58 & 654.40 & [\textbf{116}, 116] & 8.60 & 8.53 & 116 \\
30-c-c 	& 96 & 0.65 & 96 & 3.76 		& [83.78, 96] & - 		&[83.78,	96]	&	 -	&	1827.64& [\textbf{96}, 96] & 2.30 & 2.22 & [\textbf{96}, 96] & 3.22 & 3.14 & 96 \\
\hline
\end{tabular}
\begin{tablenotes}
\item[a] CPU AMD Ryzen 3700X - 4x3.6 GHz, 4x4.4 GHz, 16 threads; RAM 32 GB; {best results over  30 runs}
\item[b] CPU AMD Ryzen 3700X - 4x3.6 GHz, 4x4.4 GHz, 16 threads; RAM 32 GB; CPLEX 12.1; 3600 sec {time limit}
\item[c] CPU Intel Core i7 12700F - 4x3.6 GHz, 8x4.9 GHz, 20 threads; RAM 32 GB; Gurobi 10.0; 3600 sec  {time limit}
\item[d] CPU Intel Core i7 12700F - 4x3.6 GHz, 8x4.9 GHz, 20 threads; RAM 32 GB; OR-Tools CP-SAT 9.6; 3600 sec  {time limit}
\end{tablenotes}
\end{threeparttable}	
}
\end{table}
\end{center}}
\end{landscape}

From the results displayed in {Table \ref{tspts}, we see that inequality (\ref{va}) is  very effective in improving {the performance of the models, both $CP1$ and $MILP$.  Given this evidence, we will always consider inequalities (\ref{va}) for the next experiments.}

Table \ref{tspts}  reveals that the CP-based approach matches (or improves in the case of instance 24-rc-c) all the best-known heuristic solutions and outperforms the {exact MILP method}, both in terms of quality and times. Additionally, we observe  how the $CP1$+(\ref{va}) improves several lower bounds and closes all but three instances. 

{
To better highlight the differences between the MILP and the CP methods, we report in Figure \ref{f2} their percentage optimality gaps, calculated as $100 \cdot \frac{UB - LB}{UB}$, and in Figure \ref{f3} their required time to find the best solution (Sec$_{\text{bst}}$).
Figure \ref{f2} shows that the $CP1$ model clearly leads to lower optimality gaps than the $MILP$ model with a time limit of 3600 seconds, the latter also demonstrating scalability issues on larger instances as remarked by its linearly increasing trend (dashed line).
Whereas Figure \ref{f3} shows that the $CP1$ model is substantially faster in retrieving the best heuristic solution. These results suggest that the Constraint Programming-based approach has great potential for the PDSTSP-c problem.}

\begin{figure}[H]
\begin{center}
    \includegraphics[width=.7\linewidth]{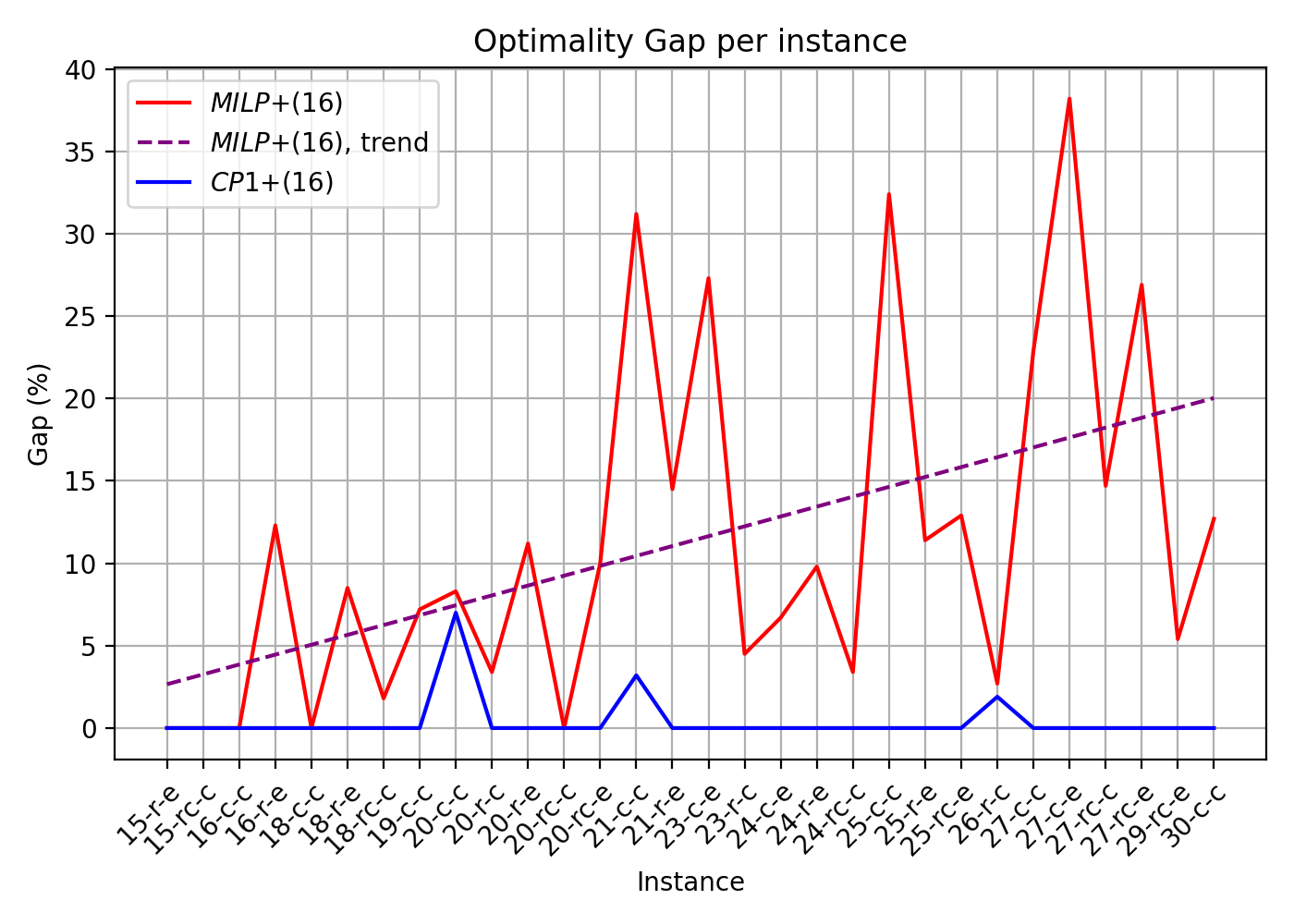}
    \caption{The optimality gap in percentage for the $MILP$+(\ref{va}) and $CP1$+(\ref{va}) on small PDSTSP-c instances.} \label{f2}
\end{center}
\end{figure}
\begin{figure}[H]
\begin{center}
\includegraphics[width=.7\linewidth]{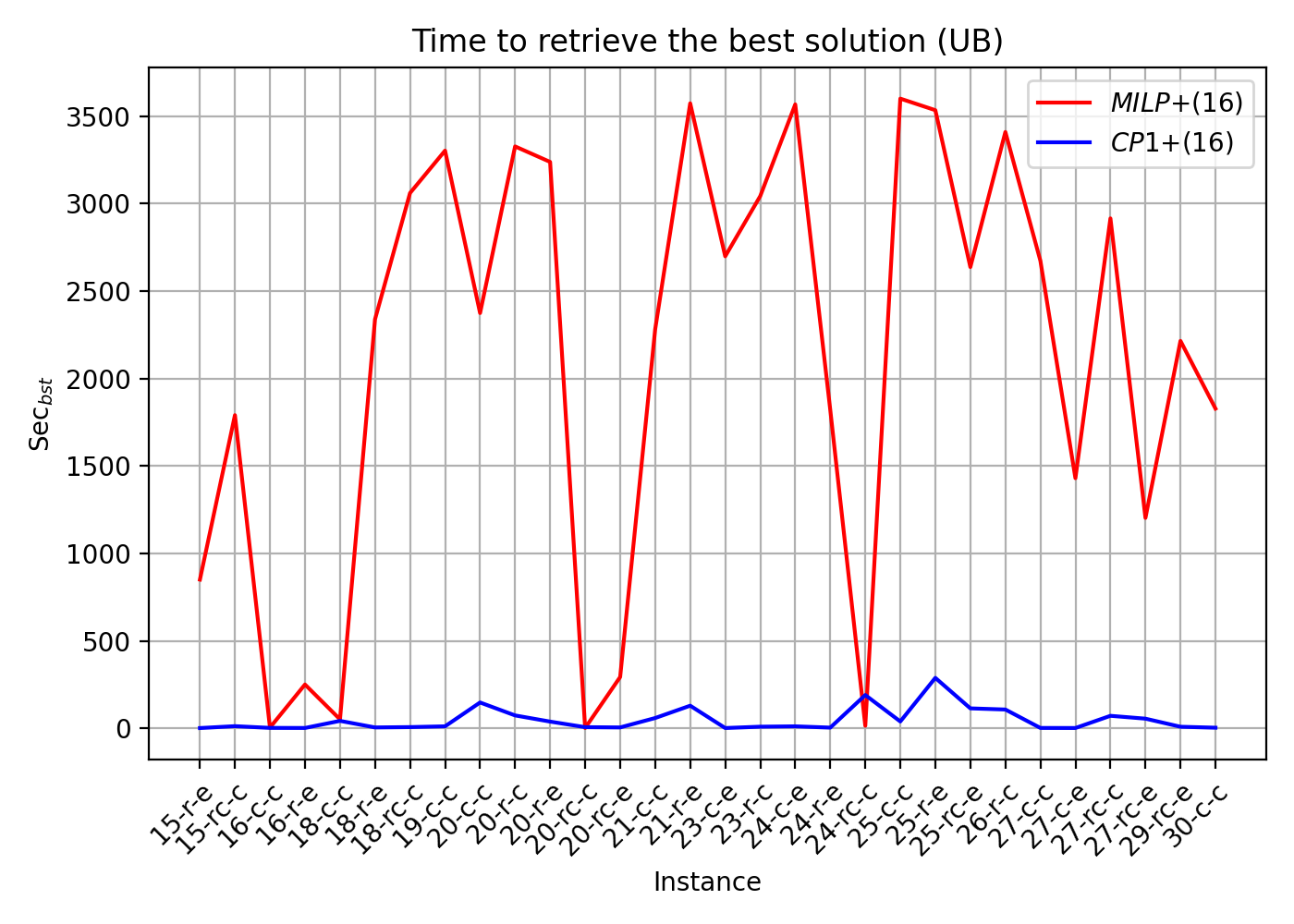}
  \caption{The time required in seconds by the $MILP$+(\ref{va})  and $CP1$+(\ref{va}) to retrieve the best heuristic solution (UB).} \label{f3}
\end{center}
\end{figure}

\begin{center}
\begin{table}[H]
\caption{Experimental results on the PDSTSP-c. Large instances.}\label{tsptl}
\resizebox{\textwidth}{!}{
\begin{threeparttable}
\begin{tabular}{l  cr  cr  crr  crr  c}
\hline
& \multicolumn{2}{c}{RnR fast  \cite{viet}\tnote{a}} & \multicolumn{2}{c}{RnR  \cite{viet}\tnote{a}} & \multicolumn{3}{c}{$CP1$\tnote{b}} & \multicolumn{3}{c}{$CP1$+(\ref{va})\tnote{b}} & Best \\
Instance & UB & Sec$_{\text{bst}}$ & UB & Sec$_{\text{bst}}$ & [LB, UB] & Sec$_{\text{tot}}$ & Sec$_{\text{bst}}$  & [LB, UB] & Sec$_{\text{tot}}$ & Sec$_{\text{bst}}$  & bounds\\
\cmidrule(lr){1-1}\cmidrule(lr){2-3}\cmidrule(lr){4-5}\cmidrule(lr){6-8}\cmidrule(lr){9-11}\cmidrule(lr){12-12}
50-r-e &  120 & 2.45 & 116 & 20.00 & [49, \textit{140}] & - & 2123.27 & [\textbf{116 }, \textit{128}] & - & 2083.72 & 116 \\
53-r-e &  136 & 2.23 & 132 & 15.17 & [68, \textit{160}] & - & 2631.50 & [\textbf{132}, \textit{140}] & - & 3177.59 & 132 \\
66-rc-e &  128 & 3.13 & 124 & 16.61 & [80, \textit{168}] & - & 393.29 & [\textbf{124}, \textit{132}] & - & 2927.03 & 124 \\
67-c-c &   76 & 2.32 & 76 & 29.05 & [68, \textit{80}] & - & 133.82 & [\textbf{73}, \textit{80}] & - & 403.51 & [73, 76] \\
68-rc-c &  79 & 2.12 & 76 & 20.45 & [72, \textit{112}] & - & 2488.01 & [\textbf{76}, 76] & 3089.06 & 3088.37 & 76 \\
76-c-c &  52 & 2.15 & 52 & 10.46 & [40, 52] & - & 26.87 & [\textbf{52}, 52] & 379.86 & 31.99 & 52 \\
82-c-e &  64 & 2.63 & 64 & 18.07 & [\textbf{64}, 64] & 33.30 & 30.55 & [\textbf{64}, 64] & 59.21 & 58.42 & 64 \\
82-rc-c &  108 & 2.65 & 104 & 27.42 & [100, \textit{140}] & - & 2373.65 & [\textbf{104}, \textit{144}] & - & 3269.15 & 104 \\
88-c-e &  108 & 3.56 & 108 & 9.62 & [\textbf{108}, 108] & 51.03 & 50.09 & [\textbf{108}, 108]  & 92.47 & 91.56 & 108 \\
91-r-c &  128 & 3.64 & 124 & 72.49 & [116, \textit{188}] & - & 2887.23 & [\textbf{120}, \textit{164}] & - & 3302.76 & [120, 124] \\
99-rc-c & 108 & 3.64 & 100 & 51.21 & [80, \textit{168}] & - & 1640.06 & [\textbf{98}, \textit{164}] & - & 2575.79 & [98, 100] \\
101-r-c &  124 & 3.98 & 120 & 99.71 & [96, \textit{176}] & - & 3456.53 & [\textbf{114}, \textit{180}] & - & 3338.22 & [114, 120] \\
103-rc-c &  128 & 4.38 & 124 & 94.25 & [108, \textit{176}] & - & 1967.67 & [\textbf{120}, \textit{164}] & - & 2147.63 & [120, 124] \\
105-rc-e &  124 & 5.78 & 120 & 62.71 & [80, \textit{184}] & - & 903.40 & [\textbf{109}, \textit{132}] & - & 1831.25 & [109, 120] \\
108-rc-e &  144 & 5.56 & 136 & 61.71 & [112, \textit{188}] & - & 3088.15 & [\textbf{134}, \textit{188}] & - & 2521.07 & [134, 136] \\
114-rc-c &  100 & 4.69 & 96 & 62.73 & [68, \textit{152}] & - & 2591.56 & [\textbf{94}, \textit{148}] & - & 2821.57 & [94, 96] \\
121-rc-e &  128 & 6.41 & 124 & 68.74 & [108, \textit{180}] & - & 2865.31 & [\textbf{121}, \textit{160}] & - & 1201.60& [121, 124] \\
126-r-c &  161 & 5.89 & 160 & 120.93 & [104, \textit{228}] & - & 1457.72 & [\textbf{151}, \textit{216}] & - & 3028.00 & [151, 160] \\
126-rc-e & 148 & 7.38 & 144 & 71.04 & [124, \textit{196}] & - & 2945.63 & [\textbf{136}, \textit{188}] & - & 1053.63 & [136, 144] \\
144-rc-r &  132 & 6.52 & 128 & 172.75 & [120, \textit{188}] & - & 2801.00 & [\textbf{122}, \textit{204}] & - & 3289.31 & [122, 128] \\
154-c-c &  72 & 7.17 & 72 & 62.96 & [\textbf{68}, 72] & - & 63.64 & [\textbf{68}, 72] & - & 61.68 & [68, 72] \\
165-r-c &  176 & 7.61 & 164 & 280.18 & [118, \textit{ 292}] & - & 2386.33 & [\textbf{140}, \textit{312}] & - & 3528.39 & [140, 164] \\
167-r-e &  200 & 10.55 & 188 & 228.31 & [72, \textit{296}] & - & 3309.32 & [\textbf{160}, \textit{304}] & - & 2209.10 & [160, 188] \\
173-r-c &  180 & 8.64 & 164 & 373.43 & [92, \textit{312}] & - & 2439.90 & [\textbf{141}, \textit{280}] & - & 2184.87 & [141, 164] \\
173-rc-r &  144 & 9.22 & 133 & 135.20 & [51, \textit{208}] & - & 2539.43 & [\textbf{115}, \textit{208}] & - & 3556.41 & [115, 133] \\
181-r-e &  232 & 11.20 & 224 & 196.09 & [125, \textit{332}] & - & 3252.35 & [\textbf{199}, \textit{348}] & - & 3566.96 & [199, 224] \\
185-c-c & 96 & 11.26 & 96 & 61.53 & [\textbf{96}, 96] & 1279.96 & 1276.68 & [\textbf{96}, 96] & 622.29 & 619.07 & 96 \\
187-rc-e &  200 & 12.67 & 196 & 119.95 & [78, \textit{ 284}] & - & 2464.47 & [\textbf{167}, \textit{288}] & - & 3228.65 & [167, 196] \\
198-c-c & 64 & 11.38 & 64 & 94.12 & [\textbf{64}, \textit{68}] & - & 82.40 & [\textbf{64}, \textit{68}] & - & 155.40 & 64 \\
200-r-e &  224 & 13.88 & 212 & 368.97 & [40, \textit{324}] & - & 3541.42 & [\textbf{162}, \textit{328}] & - & 2132.93 & [162, 212] \\
\hline
\end{tabular}
\begin{tablenotes}
\item[a] CPU AMD Ryzen 3700X - 4x3.6 GHz, 4x4.4 GHz, 16 threads; RAM 32 GB; {best results over  30 runs}
\item[b] CPU Intel Core i7 12700F - 4x3.6 GHz, 8x4.9 GHz, 20 threads; RAM 32 GB; OR-Tools CP-SAT 9.6; 3600 sec  {time limit}
\end{tablenotes}
\end{threeparttable}	
}
\end{table}
\end{center}
{
Moving to the larger instances reported in Table \ref{tsptl} we observe that $CP1$+(\ref{va}) is able to provide, for the first time, valid lower bounds for all instances. Moreover 10 over 30 bounds equal to the best known solution, hence proving for the first time the optimality of these  solutions. In the remaining instances the gaps, between the lower bound and the heuristic solution is generally small. However the upper bound provided by the CP models is not competitive with respect to that of the metaheuristic methods. Also the running times are larger, although it is worth to observe once again that for the \emph{RnR} methods the best results over 30 runs is provided, making the timing presented less fair.}

\subsection{PDSVRP-c} \label{evrp}

In this section, we aim at comparing the performance of the models $CP2$, $CP3$ and $MILP$ described in Sections~\ref{vrp2}-\ref{milp_vrp3} for the PDSVRP-c. Their results are summarized in Tables \ref{tab:vrp_small2} and \ref{tab:vrp_small3} for the small instances, covering respectively 2 and 3 trucks, and in Tables \ref{tab:vrp_large2}-\ref{tab:vrp_large5} for the large instances, using respectively 2, 3, 4, and 5 trucks. The PDSVRP-c is first introduced in this paper, so no comparison is available with  methods from other authors.

\begin{table}[h!]
\caption{Experimental results on the PDSVRP-c. Small instances, 2 trucks.}\label{tab:vrp_small2}
\resizebox{\textwidth}{!}{
\begin{threeparttable}{
\begin{tabular}{l crr crr crr c}
    \hline
    &  \multicolumn{3}{c}{$MILP$+(\ref{va})\tnote{a}}
    & \multicolumn{3}{c}{$CP2$+(\ref{va})\tnote{b}} & \multicolumn{3}{c}{$CP3$+(\ref{va})\tnote{b}} & Best \\
    Instance & [LB, UB] & Sec$_{\text{tot}}$  & Sec$_{\text{bst}}$
    & [LB, UB] & Sec$_{\text{tot}}$  & Sec$_{\text{bst}}$
    & [LB, UB] & Sec$_{\text{tot}}$  & Sec$_{\text{bst}}$  & bounds \\
    \cmidrule(lr){1-1}\cmidrule(lr){2-4}\cmidrule(lr){5-7}\cmidrule(lr){8-10}\cmidrule(lr){11-11}
    \multirow{1}{*}{15-r-e} &  		[\textbf{92	},	\textbf{92}]	&	2654.11	&	745.82 		& [\textbf{92}, \textbf{92}] & 157.1 & 0.21 & [\textbf{92}, \textbf{92}] & 2.49 & 0.45 & 92 \\
    \multirow{1}{*}{15-rc-c} & 	[\textbf{33},	\textbf{33}]	&	9.42	&	7.46 					& [\textbf{33}, \textbf{33}] & 0.92 & 0.52 & [\textbf{33}, \textbf{33}] & 5.14 & 2.92 & 33 \\
    \multirow{1}{*}{16-c-c} &  		[\textbf{40},	\textbf{40}]	&	50.04	&	3.42				& [\textbf{40}, \textbf{40}] & 6.43 & 0.29 & [\textbf{40}, \textbf{40}] & 7.17 & 2.26 & 40 \\
    \multirow{1}{*}{16-r-e} &  		[104,	\textbf{108}]	&	 -	&	783.12 				& [\textbf{108}, \textbf{108}] & 83.49 & 0.69 & [\textbf{108}, \textbf{108}] & 5.03 & 1.75 & 108 \\
    \multirow{1}{*}{18-c-c} &  		[\textbf{44},	\textbf{44}]	&	22.76	&	22.72					& [\textbf{44}, \textbf{44}] & 3.09 & 2.4 & [\textbf{44}, \textbf{44}] & 7.45 & 5.98 & 44 \\
    \multirow{1}{*}{18-r-e} &  		[\textbf{92},	\textbf{92}]	&	1153.19	&	197.42 				& [\textbf{92}, \textbf{92}] & 43.36 & 1.57 & [\textbf{92}, \textbf{92}] & 6.92 & 4.52 & 92 \\
    \multirow{1}{*}{18-rc-c} & 	[44,	\textbf{46}]	&	 -	&	3569.88  							& [\textbf{46}, \textbf{46}] & 459.47 & 161.29 & [\textbf{46}, \textbf{46}] & 401.32 & 45.46 & 46 \\
    \multirow{1}{*}{19-c-c} &  		[34, \textbf{36}]	&	 -	&	3220.87							& [\textbf{36}, \textbf{36}] & 7.85 & 4.15 & [\textbf{36}, \textbf{36}] & 15.41 & 3.84 & 36 \\
    \multirow{1}{*}{20-c-c} &  		[\textbf{40},	\textbf{40}]	&	907.04	&	5.25					& [\textbf{40}, \textbf{40}] & 4.35 & 2.23 & [\textbf{40}, \textbf{40}] & 4.26 & 1.75 & 40 \\
    \multirow{1}{*}{20-r-c} &  		[\textbf{48},	\textbf{48}]	&	2023.70	&	1886.70  							& [\textbf{48}, \textbf{48}] & 454.46 & 9.19 & [\textbf{48}, \textbf{48}] & 62.72 & 43.56 & 48 \\
    \multirow{1}{*}{20-r-e} &  		[63,		76]	&	 -	&	2707.00 					& [\textbf{72}, \textbf{72}] & 462.23 & 331.95 & [\textbf{72}, \textbf{72}] & 27.42 & 23.27 & 72 \\
    \multirow{1}{*}{20-rc-c} & 	[63,	\textbf{64}]	&	 -	&	1977.44  					& [58, \textbf{64}] & - & 1.15 & [\textbf{64}, \textbf{64}] & 14.57 & 8.09 & 64 \\
    \multirow{1}{*}{20-rc-e} & 	[72,	\textbf{80}]&	 -	&	3031.82 					& [\textbf{64}, \textbf{80}] & - & 1.49 & [\textbf{80}, \textbf{80}] & 24.93 & 7.82 & 80 \\
    \multirow{1}{*}{21-c-c} &  		[\textbf{40},	\textbf{40}]	&	59.12	&	8.87  					& [\textbf{40}, \textbf{40}] & 9.94 & 1.36 & [\textbf{40}, \textbf{40}] & 11.6 & 2.01 &  40 \\
    \multirow{1}{*}{21-r-e} &  		[51,	\textbf{76}]	&	 -	&	3348.26 							& [\textbf{76}, \textbf{76}] & 475.63 & 4.08 & [\textbf{76}, \textbf{76}] & 82.46 & 9.41 & 76 \\
    \multirow{1}{*}{23-c-e} & 		[44,	\textbf{80}]	&	 -	&	0.98					& [42, \textbf{80}] & - & 0.98 & [\textbf{80}, \textbf{80}] & 17.85 & 0.67 & 80 \\
    \multirow{1}{*}{23-r-c} &  		[\textbf{60},	\textbf{60}]	&	1536.90	&	1379.17 							& [57, \textbf{60}] & - & 9.64 & [\textbf{60}, \textbf{60}] & 475.01 & 10.69 & 60 \\
    \multirow{1}{*}{24-c-e} & 		[56,	\textbf{60}]	&	 -	&	505.00 					& [\textbf{60}, \textbf{60}] & 116.75 & 34.53 & [\textbf{60}, \textbf{60}] & 54.04 & 53.77 & 60 \\
    \multirow{1}{*}{24-r-e} &  		[68,	\textbf{100}]	&	 -	&	0.89  							& [67, \textbf{100}] & - & 3.09 & [\textbf{100}, \textbf{100}] & 68.82 & 35.52 & 100 \\
    \multirow{1}{*}{24-rc-c} & 	[49,	\textbf{52}]	&	 -	&	42.34  							& [\textbf{52}, \textbf{52}] & 1350.35 & 25.78 & [\textbf{52}, \textbf{52}] & 307.8 & 82.07 & 52 \\
    \multirow{1}{*}{25-c-c} &  		[\textbf{40},	\textbf{40}]	&	2523.46	&	1988.78 							& [\textbf{40}, \textbf{40}] & 899.8 & 26.8 & [\textbf{40}, \textbf{40}] & 79.98 & 63.91 & 40 \\
    \multirow{1}{*}{25-r-e} &  		[72,		92]	&	 -	&	693.55 					& [85, \textbf{88}] & - & 35.41 & [\textbf{88}, \textbf{88}] & 175.75 & 44.11 & 88 \\
    \multirow{1}{*}{25-rc-e} & 	[64, 80]	&	 -	&	239.90 					& [59, \textbf{76}] & - & 80.16 & [\textbf{76}, \textbf{76}] & 189.81 & 3.42 & 76 \\
    \multirow{1}{*}{26-r-c} &   		[68,	\textbf{70}]	&	 -	&	1381.35 					& [65, \textbf{70}] & - & 20.87 & [\textbf{70}, \textbf{70}] & 2701.81 & 700.35 & 70 \\
    \multirow{1}{*}{27-c-c} & 		[40,	\textbf{52}]	&	 -	&	83.87				& [42, \textbf{52}] & - & 2.32 & [\textbf{52}, \textbf{52}] & 37.07 & 16.67 & 52 \\
    \multirow{1}{*}{27-c-e} & 		[7,	\textbf{68}]	&	-	&	6.18 					& [\textbf{68}, \textbf{68}] & 2741.9 & 0.65 & [\textbf{68}, \textbf{68}] & 130.71 & 1.56 & 68 \\
    \multirow{1}{*}{27-rc-c} & 	[64,	\textbf{72}]&	 -	&	58.13 					& [\textbf{64}, \textbf{72}] & - & 8.96 & [\textbf{72}, \textbf{72}] & 79.95 & 53.96 &  72 \\
    \multirow{1}{*}{27-rc-e} & 	[36,		80]	&	 -	&	171.86					& [47, \textbf{76}] & - & 20.22 & [\textbf{76}, \textbf{76}] & 198.05 & 46.71 & 76 \\
    \multirow{1}{*}{29-rc-e} & 	[62,	\textbf{100}]	&	 -	&	3092.20							& [65, \textbf{100}] & - & 13.71 & [\textbf{100}, \textbf{100}] & 59.26 & 48.57 & 100 \\
    \multirow{1}{*}{30-c-c} &  		[36,	\textbf{64}]	&	 -	&	2935.06 					& [48, \textbf{64}] & - & 2.09 & [\textbf{64}, \textbf{64}] & 39.64 & 2.99 & 64 \\  \hline
\end{tabular}}
\begin{tablenotes}
\item[a] CPU Intel Core i7 12700F - 4x3.6 GHz, 8x4.9 GHz, 20 threads; RAM 32 GB; Gurobi 10.0; 3600 sec  {time limit}
\item[b] CPU Intel Core i7 12700F - 4x3.6 GHz, 8x4.9 GHz, 20 threads; RAM 32 GB; OR-Tools CP-SAT 9.6; 3600 sec  {time limit}
\end{tablenotes}
\end{threeparttable}	
}
\end{table}
\begin{table}[h!]
\caption{Experimental results on the PDSVRP-c. Small instances, 3 trucks.}\label{tab:vrp_small3}
\resizebox{\textwidth}{!}{
\begin{threeparttable}{
\begin{tabular}{l crr crr crr c}
    \hline
    &  \multicolumn{3}{c}{$MILP$+(\ref{va})\tnote{a}}
    & \multicolumn{3}{c}{$CP2$+(\ref{va})\tnote{b}} & \multicolumn{3}{c}{$CP3$+(\ref{va})\tnote{b}} & Best \\
    Instance  & [LB, UB] & Sec$_{\text{tot}}$  & Sec$_{\text{bst}}$
    & [LB, UB] & Sec$_{\text{tot}}$  & Sec$_{\text{bst}}$
    & [LB, UB] & Sec$_{\text{tot}}$  & Sec$_{\text{bst}}$  & bounds \\
    \cmidrule(lr){1-1}\cmidrule(lr){2-4}\cmidrule(lr){5-7}\cmidrule(lr){8-10}\cmidrule(lr){11-11}
    \multirow{1}{*}{15-r-e} & [72,	\textbf{92}]	&	 -	&	3401.20				& [\textbf{92}, \textbf{92}] & 331.05 & 0.19 & [\textbf{92}, \textbf{92}] & 4.44 & 0.36 & 92 \\
    \multirow{1}{*}{15-rc-c} & [\textbf{32}, \textbf{32}]	&	11.08	&	11.01				& [\textbf{32}, \textbf{32}] & 1.43 & 1.27 & [\textbf{32}, \textbf{32}] & 5.95 & 3.66 & 32 \\
    \multirow{1}{*}{16-c-c} & [\textbf{36},	\textbf{36}]	&	2.20	&	2.18 		& [\textbf{36}, \textbf{36}] & 1.21 & 1.17 & [\textbf{36}, \textbf{36}] & 6.01 & 5.41 &  36\\
    \multirow{1}{*}{16-r-e} & [68,	\textbf{108}]	&	 -	&	2332.97				& [\textbf{108}, \textbf{108}] & 290.7 & 1.05 & [\textbf{108}, \textbf{108}] & 6.55 & 1.15 & 108 \\
    \multirow{1}{*}{18-c-c} & [\textbf{44},	\textbf{44}]	&	25.27	&	25.21 						& [\textbf{44}, \textbf{44}] & 3.02 & 1.61 & [\textbf{44}, \textbf{44}] & 6.37 & 2.12 & 44 \\
    \multirow{1}{*}{18-r-e} & [88,	\textbf{92}]	&	 -	&	614.09 				& [\textbf{92}, \textbf{92}] & 22.3 & 0.92 & [\textbf{92}, \textbf{92}] & 12 & 2.1 & 92 \\
    \multirow{1}{*}{18-rc-c} & [\textbf{40},	\textbf{40}]	&	18.60	&	18.53  								& [\textbf{40}, \textbf{40}] & 5.22 & 5.09 & [\textbf{40}, \textbf{40}] & 30.65 & 24.69 & 40 \\
    \multirow{1}{*}{19-c-c} & [29,	\textbf{36}]	&	 -	&	3554.00  								& [\textbf{36}, \textbf{36}] & 2.26 & 0.86 & [\textbf{36}, \textbf{36}] & 19.25 & 3.59 & 36 \\
    \multirow{1}{*}{20-c-c} & [\textbf{40},	\textbf{40}]	&	383.19	&	139.78  								& [\textbf{40}, \textbf{40}] & 1.9 & 0.72 & [\textbf{40}, \textbf{40}] & 7.35 & 1.27 & 40 \\
    \multirow{1}{*}{20-r-c} & [\textbf{37},	\textbf{37}]	&	1202.72	&	582.24  								& [\textbf{37}, \textbf{37}] & 8.3 & 3.81 & [\textbf{37}, \textbf{37}] & 99.17 & 36.74 & 37 \\
    \multirow{1}{*}{20-r-e} & [44,	\textbf{72}]	&	 -	&	45.40 								& [\textbf{72}, \textbf{72}] & 443.37 & 7.4 & [\textbf{72}, \textbf{72}] & 17.12 & 9.27 & 72 \\
    \multirow{1}{*}{20-rc-c} & [\textbf{48},	\textbf{48}]	&	604.42	&	8.91  								& [\textbf{48}, \textbf{48}] & 99.53 & 3.15 & [\textbf{48}, \textbf{48}] & 94.88 & 26.53 & 48 \\
    \multirow{1}{*}{20-rc-e} & [60,	\textbf{68}]	&	 -	&	595.63  								& [\textbf{68}, \textbf{68}] & 122.18 & 42 & [\textbf{68}, \textbf{68}] & 45.74 & 8.25 & 68 \\
    \multirow{1}{*}{21-c-c} & [\textbf{36},	\textbf{36}]	&	28.87	&	6.92 								& [\textbf{36}, \textbf{36}] & 4.96 & 4.75 & [\textbf{36}, \textbf{36}] & 16.5 & 7.81 & 36 \\
    \multirow{1}{*}{21-r-e} & [34,	\textbf{76}]	&	 -	&	689.44 								& [\textbf{76}, \textbf{76}] & 891.56 & 11.3 & [\textbf{76}, \textbf{76}] & 427.99 & 28.05 & 76 \\
    \multirow{1}{*}{23-c-e} & [36,	\textbf{80}]	&	 -	&	3.88  								& [66, \textbf{80}] & - & 0.94 & [\textbf{80}, \textbf{80}] & 25.44 & 0.93 & 80 \\
    \multirow{1}{*}{23-r-c} & [\textbf{48},	\textbf{48}]	&	203.22	&	44.94 								& [\textbf{48}, \textbf{48}] & 20.92 & 14.24 & [\textbf{48}, \textbf{48}] & 590.1 & 135.03 & 48 \\
    \multirow{1}{*}{24-c-e} & [56,	\textbf{60}]	&	 -	&	2000.96  								& [\textbf{60}, \textbf{60}] & 135.61 & 3.35 & [\textbf{60}, \textbf{60}] & 28.56 & 11.45 & 60 \\
    \multirow{1}{*}{24-r-e} & [47,	\textbf{100}]	&	 -	&	8.28  								& [80, \textbf{100}] & - & 1.72 & [\textbf{100}, \textbf{100}] & 64.88 & 8.77 & 100 \\
    \multirow{1}{*}{24-rc-c} & [41,	\textbf{44}]	&	 -	&	3456.94  								& [\textbf{44}, \textbf{44}] & 227.22 & 77.94 & [\textbf{44}, \textbf{44}] & 926.73 & 387.43 & 44 \\
    \multirow{1}{*}{25-c-c} & [30,		40]	&	 -	&	97.39 								& [\textbf{37}, \textbf{37}] & 70.58 & 31.83 & [\textbf{37}, \textbf{37}] & 107.3 & 42.94 & 37 \\
    \multirow{1}{*}{25-r-e} & [57,	96]	&	 -	&	3298.63					& [59, \textbf{85}] & - & 44.6 & [\textbf{85}, \textbf{85}] & 404.43 & 242.54 & 85 \\
    \multirow{1}{*}{25-rc-e} & [52,		69]	&	 -	&	2601.65					& [65, \textbf{66}] & - & 141.25 & [\textbf{66}, \textbf{66}] & 356.83 & 61.26 & 66 \\
    \multirow{1}{*}{26-r-c} & [\textbf{56},	\textbf{56}]	&	180.88	&	67.14					& [52, \textbf{56}] & - & 1399.21 & [{55}, \textbf{56}] & - & 70.15 &  56 \\
    \multirow{1}{*}{27-c-c} & [\textbf{36},	\textbf{36}]	&	1367.86	&	936.48					& [\textbf{36}, \textbf{36}] & 577.74 & 3.69 & [\textbf{36}, \textbf{36}] & 121.7 & 39.48 & 36 \\
    \multirow{1}{*}{27-c-e} & [8,	\textbf{68}]	&	 -	&	4.92					& [\textbf{68}, \textbf{68}] & 2315.01 & 1.31 & [\textbf{68}, \textbf{68}] & 437.75 & 1.71 & 68 \\
    \multirow{1}{*}{27-rc-c} & [\textbf{60},	\textbf{60}]	&	223.20	&	213.85 					& [\textbf{60}, \textbf{60}] & 6.79 & 6.01 & [\textbf{60}, \textbf{60}] & 74.79 & 63.9 & 60 \\
    \multirow{1}{*}{27-rc-e} & [28,	\textbf{76}]	&	 -	&	155.68 					& [56, \textbf{76}] & - & 8.02 & [\textbf{76}, \textbf{76}] & 520.4 & 14.61 & 76 \\
    \multirow{1}{*}{29-rc-e} & [53,		108]	&	 -	&	3337.59							& [72, \textbf{100}] & - & 23.67 & [\textbf{100}, \textbf{100}] & 56.98 & 35.42 & 100 \\
    \multirow{1}{*}{30-c-c} & [26,	\textbf{38}	]	&	 -	&	3182.76 					& [\textbf{38}, \textbf{38}] & 96.12 & 5.63 & [\textbf{38}, \textbf{38}] & 145.11 & 15.99 & 38 \\
 \hline
\end{tabular}}
\begin{tablenotes}
\item[a] CPU Intel Core i7 12700F - 4x3.6 GHz, 8x4.9 GHz, 20 threads; RAM 32 GB; Gurobi 10.0; 3600 sec  {time limit}
\item[b] CPU Intel Core i7 12700F - 4x3.6 GHz, 8x4.9 GHz, 20 threads; RAM 32 GB; OR-Tools CP-SAT 9.6; 3600 sec  {time limit}
\end{tablenotes}
\end{threeparttable}	
}
\end{table}

Tables \ref{tab:vrp_small2} and \ref{tab:vrp_small3} suggest that solving the $MILP$ is less effective than solving the CP models, both in terms of bounds provided and computing time. The only remarkable exception is instance {26-r-c}, which is closed by the former but not by the latters. 

{One can also observe that the performances of $MILP$ degrade with the increasing of the instance size, much more than that of the  CP methods. After some test with larger instances (not reported here), and considering the analogous decision made for the PDSTSP-c in \cite{viet}, we decided to not consider the $MILP$ for the experiments on large instances (Tables \ref{tab:vrp_large2}-\ref{tab:vrp_large5})}. 

The results of the two CP models suggest that the 3-indices formulation ($CP3$) is superior, being able to close all the instances but one. The 2-indices model appears slower even though the quality of its upper bounds is {the same of $CP3$}. This highlights that the weakness of the $CP2$ model is {in the computation of} the lower bound. 

The results reported in Tables \ref{tab:vrp_large2}-\ref{tab:vrp_large5} for large instances (note that the column Sec$_{\text{tot}}$ has been omitted, since no optimality is proven) and a varying number of trucks lead to the following observations. The model with 3 indices, which performs the best on small instances (see Tables {\ref{tab:vrp_small2} and \ref{tab:vrp_small3}), is instead performing worse than the 2-indices model on large ones, especially in terms of retrieved lower bounds. This might suggest that handling multiple truck tours with the \emph{MultipleCircuit} command becomes effective when tours are complex.

{There are however a few exceptions where the 3-indices model is better either in terms of lower or upper bounds. Specifically, the $CP3$ model appears to be more consistent in instances with many customers and a few trucks, in which the 2-indices model often fails to produce any feasible solution. This might indicate that the models are approaching their natural limit.}

\begin{table}[H]
{
\begin{center}
\caption{Experimental results on the PDSVRP-c. Large instances, 2 trucks.}\label{tab:vrp_large2}
\resizebox{0.66\textwidth}{!}{
\begin{threeparttable}{
\begin{tabular}{l cr cr c}
\hline
& \multicolumn{2}{c}{$CP2$+(\ref{va})\tnote{c}} & \multicolumn{2}{c}{$CP3$+(\ref{va})\tnote{c}} & Best \\
Instance & [LB, UB] & Sec$_{\text{bst}}$ & [LB, UB] & Sec$_{\text{bst}}$ & bounds \\
\cmidrule(lr){1-1}\cmidrule(lr){2-3}\cmidrule(lr){4-5}\cmidrule(lr){6-6}
50-r-e  & [\textbf{65}, \textbf{116}]   & 206.57 & [63,  120]   & 168.48 & [65, 116] \\
53-r-e  & [77, \textbf{112}]   & 894.09 & [\textbf{82},  128]   & 1756.80 & [82, 112] \\
66-rc-e  & [72, \textbf{112}]   & 1829.73 & [\textbf{73},  136]   & 866.28 & [73, 112] \\
67-c-c  & [\textbf{38}, \textbf{52}]   & 22.33 & [31,  \textbf{52}]   & 827.01 & [38, 52] \\
68-rc-c  & [50, \textbf{56}]   & 3332.51 & [\textbf{52},  104]   & 3088.50 & [52, 56] \\
76-c-c  & [\textbf{26}, \textbf{36}]   & 20.60 & [16,  40]   & 185.95 & [26, 36] \\
82-c-e  & [\textbf{32}, \textbf{64}]   & 25.41 & [17,  \textbf{64}]   & 73.68 & [32, 64] \\
82-rc-c  & [\textbf{62}, \textbf{116}]   & 2974.84 & [56,  132]   & 2615.62 & [62, 116] \\
88-c-e  & [54, \textbf{84}]   & 298.18 & [\textbf{58},  112]   & 49.28 & [58, 84] \\
91-r-c  & [\textbf{75}, \textbf{152}]   & 405.02 & [\textbf{75},  160]   & 2249.67 & [75, 152] \\
99-rc-c  & [\textbf{63}, \textbf{96}]   & 2083.65 & [51,  144]   & 564.95 & [63, 96] \\
101-rc  & [\textbf{71}, 164]   & 2921.49 & [53,  \textbf{152}]   & 1731.45 & [71, 152] \\
103-rc-c  & [\textbf{69}, \textbf{124}]   & 2603.95 & [52,  128]   & 2912.93 & [69, 124] \\
105-rc-e  & [\textbf{65}, \textbf{136}]   & 2170.84 & [57,  148]   & 1383.74 & [65, 136] \\
108-rc-e  & [\textbf{79}, 172]   & 1683.07 & [70,  \textbf{160}]   & 831.13 & [79, 160] \\
114-rc-c  & [\textbf{58}, \textbf{124}]   & 3417.23 & [49,  140]   & 411.62 & [58, 124] \\
121-rc-e  & [\textbf{70}, 156]   & 647.12 & [56,  \textbf{152}]   & 2088.27 & [70, 152] \\
126-rc-e  & [\textbf{87}, 220]   & 3115.59 & [67,  \textbf{184}]   & 1956.96 & [87, 184] \\
126-r-c  & [\textbf{78}, 160]   & 2679.11 & [56,  \textbf{156}]   & 1448.65 & [78, 156] \\
144-rc-c  & [\textbf{67}, 272]   & 2610.83 & [47,  \textbf{168}]   & 3103.46 & [67, 168] \\
154-c-c  & [\textbf{35}, -]  & -   & [8,  \textbf{72}]   & 279.16 & [35, 72] \\
165-r-c  & [\textbf{88}, -]  & -   & [67,  \textbf{224}]  & 3544.74 & [88, 224] \\
167-r-e  & [\textbf{100}, -] & - & [74,  \textbf{256}]  & 3151.22 & [100, 256] \\
173-r-c  & [\textbf{85}, \textbf{204}]  & 2929.34 & [59,  240]  & 2251.40 & [85, 204] \\
173-rc-c  & [\textbf{79}, -] & - & [48,  \textbf{180}]  & 1797.98 & [79, 180] \\
181-r-e  & [\textbf{112}, -] & - & [78,  \textbf{252}]  & 3388.32 & [112, 252] \\
185-c-c  & [\textbf{48}, -] & - & [24,  \textbf{96}]  & 316.31 & [48, 96] \\
187-rc-e  & [\textbf{100}, 308]  & 3391.71 & [65,  \textbf{212}]  & 1567.38 & [100, 212] \\
198-c-c  & [\textbf{32}, -] & - & [12,  \textbf{64}]  & 271.52 & [32, 64] \\
200-r-e  & [\textbf{105}, -] & - & [68,  \textbf{324}]  & 2072.94 & [105, 324] \\
 \hline
\end{tabular}}
\begin{tablenotes}
\item[c] CPU Intel Core i7 12700F - 4x3.6 GHz, 8x4.9 GHz, 20 threads; RAM 32 GB; OR-Tools CP-SAT 9.6; 3600 sec  {time limit}
\end{tablenotes}
\end{threeparttable}	
}
\end{center}}
\end{table}
\begin{table}[H]
{
\begin{center}
\caption{Experimental results on the PDSVRP-c. Large instances, 3 trucks.}\label{tab:vrp_large3}
\resizebox{0.66\textwidth}{!}{
\begin{threeparttable}{
\begin{tabular}{l cr cr c}
\hline
& \multicolumn{2}{c}{$CP2$+(\ref{va})\tnote{c}} & \multicolumn{2}{c}{$CP3$+(\ref{va})\tnote{c}} & Best \\
Instance  & [LB, UB] & Sec$_{\text{bst}}$ & [LB, UB] & Sec$_{\text{bst}}$ & bounds \\
\cmidrule(lr){1-1}\cmidrule(lr){2-3}\cmidrule(lr){4-5}\cmidrule(lr){6-6}
50-r-e  & [\textbf{48}, \textbf{112}]   & 79.65 & [47,  112]   & 411.11 & [48, 112] \\
53-r-e  & [\textbf{56}, \textbf{96}]   & 860.00 & [51,  112]   & 2074.85 & [56, 96] \\
66-rc-e  & [\textbf{53}, \textbf{108}]   & 282.64 & [38,  116]   & 139.49 & [53, 108] \\
67-c-c  & [\textbf{27}, \textbf{52}]   & 32.82 & [9,  \textbf{52}]   & 353.35 & [27, 52] \\
68-rc-c  & [\textbf{39}, \textbf{56}]   & 756.18 & [34,  104]   & 655.52 & [39, 56] \\
76-c-c  & [\textbf{18}, \textbf{24}]   & 42.28 & [12,  52]   & 81.65 & [18, 24] \\
82-c-e  & [\textbf{22}, \textbf{64}]   & 21.88 & [8,  \textbf{64}]   & 26.79 & [22, 64] \\
82-rc-c  & [\textbf{47}, \textbf{80}]   & 1727.16 & [38,  128]   & 312.65 & [47, 80] \\
88-c-e  & [\textbf{36}, \textbf{76}]   & 375.27 & [32,  104]   & 587.30 & [36, 76] \\
91-r-c  & [\textbf{56}, \textbf{120}]   & 3036.11 & [42,  148]   & 726.56 & [56, 120] \\
99-rc-c  & [\textbf{47}, \textbf{64}]   & 2650.32 & [29,  128]   & 196.67 & [47, 64] \\
101-rc  & [\textbf{52}, \textbf{128}]   & 2645.43 & [36,  144]   & 2520.98 & [52, 128] \\
103-rc-c  & [\textbf{49}, \textbf{96}]   & 2229.84 & [32,  136]   & 2332.89 & [49, 96] \\
105-rc-e  & [\textbf{49}, \textbf{120}]   & 877.50 & [34,  132]   & 907.44 & [49, 120] \\
108-rc-e  & [\textbf{58}, 184]   & 1969.45 & [37,  \textbf{160}]   & 1273.20 & [58, 160] \\
114-rc-c  & [\textbf{44}, \textbf{80}]   & 1676.32 & [32,  112]   & 466.45 & [44, 80] \\
121-rc-e  & [\textbf{52}, \textbf{124}]   & 2820.31 & [40,  152]   & 1701.91 & [52, 124] \\
126-rc-e  & [\textbf{63}, \textbf{136}]   & 2839.24 & [44,  164]   & 2663.93 & [63, 136] \\
126-r-c  & [\textbf{56}, \textbf{140}]   & 2191.71 & [38,  148]   & 3114.44 & [56, 140] \\
144-rc-c  & [\textbf{50}, \textbf{132}]   & 3362.32 & [35,  160]   & 2396.79 & [50, 132] \\
154-c-c  & [\textbf{24}, \textbf{36}]   & 195.44 & [8,  68]   & 1368.67 & [24, 36] \\
165-r-c  & [\textbf{68}, -] & - & [50,  \textbf{212}]  & 3120.16 & [68, 212] \\
167-r-e  & [\textbf{73}, -] & - & [54,  \textbf{204}]  & 2112.65 & [73, 204] \\
173-r-c  & [\textbf{65}, -] & - & [45,  \textbf{212}]  & 2004.93 & [65, 212] \\
173-rc-c  & [\textbf{58}, 172]  & 2994.35 & [37,  \textbf{168}]  & 2592.28 & [58, 168] \\
181-r-e  & [\textbf{82}, -] & - & [55,  \textbf{216}]  & 3342.10 & [82, 216] \\
185-c-c  & [\textbf{32}, -]  & - & [14,  \textbf{96}]  & 1280.86 & [32, 96] \\
187-rc-e  & [\textbf{74}, -] & - & [46,  \textbf{212}]  & 2849.33 & [74, 212] \\
198-c-c  & [\textbf{22}, \textbf{36}]  & 158.92 & [8,  68]  & 108.97 & [22, 36] \\
200-r-e  & [\textbf{77}, -] & - & [48,  \textbf{252}]  & 1817.10 & [77, 252] \\
 \hline
\end{tabular}}
\begin{tablenotes}
\item[c] CPU Intel Core i7 12700F - 4x3.6 GHz, 8x4.9 GHz, 20 threads; RAM 32 GB; OR-Tools CP-SAT 9.6; 3600 sec  {time limit}
\end{tablenotes}
\end{threeparttable}	
}
\end{center}}
\end{table}
\begin{table}[H]
{
\begin{center}
\caption{Experimental results on the PDSVRP-c. Large instances, 4 trucks.}\label{tab:vrp_large4}
\resizebox{0.66\textwidth}{!}{
\begin{threeparttable}{
\begin{tabular}{l cr cr c}
\hline
& \multicolumn{2}{c}{$CP2$+(\ref{va})\tnote{c}} & \multicolumn{2}{c}{$CP3$+(\ref{va})\tnote{c}} & Best \\
Instance & [LB, UB] & Sec$_{\text{bst}}$  & [LB, UB] & Sec$_{\text{bst}}$  & bounds \\
\cmidrule(lr){1-1}\cmidrule(lr){2-3}\cmidrule(lr){4-5}\cmidrule(lr){6-6}
50-r-e  & [\textbf{46}, \textbf{104}]   & 649.55 & [35,  112]  & 213.25 & [46, 104] \\
53-r-e  & [\textbf{50}, \textbf{96}]   & 1068.12 & [38,  112]   & 548.64 & [50, 96] \\
66-rc-e  & [\textbf{41}, \textbf{104}]   & 3493.30 & [34,  108]   & 1019.85 & [41, 104] \\
67-c-c  & [\textbf{21}, \textbf{48}]   & 25.68 & [8,  52]   & 1297.51 & [21, 48] \\
68-rc-c  & [\textbf{32}, \textbf{52}]   & 410.57 & [29,  88]   & 296.61 & [32, 52] \\
76-c-c  & [\textbf{14}, \textbf{24}]   & 44.33 & [12,  56]   & 24.41 & [14, 24] \\
82-c-e  & [\textbf{18}, \textbf{64}]   & 18.15 & [8,  \textbf{64}]   & 20.44 & [18, 64] \\
82-rc-c  & [\textbf{38}, \textbf{68}]   & 2275.00 & [31,  124]   & 194.26 & [38, 68] \\
88-c-e  & [28, \textbf{76}]   & 76.88 & [\textbf{32},  108]   & 1177.87 & [32, 76] \\
91-r-c  & [\textbf{45}, \textbf{96}]   & 3019.26 & [32,  156]   & 248.69 & [45, 96] \\
99-rc-c  & [\textbf{37}, \textbf{68}]   & 1058.95 & [24,  120]   & 322.23 & [37, 68] \\
101-rc  & [\textbf{42}, \textbf{76}]   & 3171.49 & [30,  144]   & 2589.25 & [42, 76] \\
103-rc-c  & [\textbf{39}, \textbf{80}]   & 1490.89 & [26,  140]   & 521.14 & [39, 80] \\
105-rc-e  & [\textbf{39}, \textbf{116}]   & 261.16 & [26,  132]   & 1691.38 & [39, 116] \\
108-rc-e  & [\textbf{46}, \textbf{124}]   & 454.82 & [28,  152]   & 3163.13 & [46, 124] \\
114-rc-c  & [\textbf{35}, \textbf{88}]   & 2564.47 & [26,  120]   & 1369.36 & [35, 88] \\
121-rc-e  & [\textbf{42}, \textbf{104}]   & 3185.34 & [29,  144]   & 410.13 & [42, 104] \\
126-rc-e  & [\textbf{50}, \textbf{132}]   & 3362.14 & [35,  164]   & 2600.11 & [50, 132] \\
126-r-c  & [\textbf{45}, \textbf{116}]   & 1094.09 & [28,  140]   & 729.26 & [45, 116] \\
144-rc-c  & [\textbf{40}, \textbf{128}]   & 3013.06 & [25,  144]   & 2451.13 & [40, 128] \\
154-c-c  & [\textbf{18}, \textbf{40}]   & 949.21 & [8,  72]   & 63.77 & [18, 40] \\
165-r-c  & [\textbf{54}, \textbf{192}]  & 243.15 & [40, \textbf{192}]  & 3124.08 & [54, 192] \\
167-r-e  & [\textbf{58}, \textbf{176}]  & 3277.00 & [42,  196]  & 1489.17 & [58, 176] \\
173-r-c  & [\textbf{54}, 352]  & 3435.45 & [36,  \textbf{192}]  & 3070.29 & [54, 192] \\
173-rc-c  & [\textbf{46}, \textbf{116}]  & 1650.91 & [29,  164]  & 3368.14 & [46, 116] \\
181-r-e  & [\textbf{65}, 268]  & 2937.67 & [42,  \textbf{208}]  & 3048.86 & [65, 208] \\
185-c-c  & [\textbf{24}, \textbf{48}]  & 2350.91 & [14,  100]  & 161.08 & [24, 48] \\
187-rc-e  & [\textbf{58}, 216]  & 2551.50 & [37,  \textbf{204}]  & 2097.96 & [58, 204] \\
198-c-c  & [\textbf{16}, -] & - & [8,  \textbf{68}]  & 122.39 & [16, 68] \\
200-r-e  & [\textbf{60}, 308]  & 3550.81 & [38,  \textbf{228}]  & 2613.35 & [60, 228] \\
 \hline
\end{tabular}}
\begin{tablenotes}
\item[c] CPU Intel Core i7 12700F - 4x3.6 GHz, 8x4.9 GHz, 20 threads; RAM 32 GB; OR-Tools CP-SAT 9.6; 3600 sec {time limit}
\end{tablenotes}
\end{threeparttable}	
}
\end{center}}
\end{table}
\begin{table}[H]
{
\begin{center}
\caption{Experimental results on the PDSVRP-c. Large instances, 5 trucks.}\label{tab:vrp_large5}
\resizebox{0.66\textwidth}{!}{
\begin{threeparttable}{
\begin{tabular}{l cr cr c}
\hline
& \multicolumn{2}{c}{$CP2$+(\ref{va})\tnote{c}} & \multicolumn{2}{c}{$CP3$+(\ref{va})\tnote{c}} & Best \\
Instance & [LB, UB] & Sec$_{\text{bst}}$ & [LB, UB] & Sec$_{\text{bst}}$ & bounds\\
\cmidrule(lr){1-1}\cmidrule(lr){2-3}\cmidrule(lr){4-5}\cmidrule(lr){6-6}
50-r-e  & [\textbf{47}, \textbf{100}]   & 227.16 & [30,  112]   & 54.55 & [47, 100] \\
53-r-e  & [\textbf{50}, \textbf{92}]   & 645.51 & [32,  112]   & 667.24 & [50, 92] \\
66-rc-e  & [\textbf{35}, \textbf{100}]   & 487.89 & [24,  120]   & 482.57 & [35, 100] \\
67-c-c  & [\textbf{18}, \textbf{52}]   & 64.31 & [8,  \textbf{52}]   & 1437.96 & [18, 52] \\
68-rc-c  & [\textbf{28}, \textbf{44}]   & 1047.89 & [23,  80]   & 1776.83 & [28, 44] \\
76-c-c  & [\textbf{12}, \textbf{24}]   & 67.04 & [\textbf{12},  40]   & 400.81 & [12, 24] \\
82-c-e  & [\textbf{15}, \textbf{64}]   & 17.43 & [6,  \textbf{64}]   & 25.43 & [15, 64] \\
82-rc-c  & [\textbf{32}, \textbf{68}]   & 592.96 & [24,  112]   & 782.49 & [32, 68] \\
88-c-e  & [23, \textbf{72}]   & 218.44 & [\textbf{32},  108]   & 109.66 & [32, 72] \\
91-r-c  & [\textbf{38}, \textbf{88}]   & 3122.27 & [28,  124]   & 3272.57 & [38, 88] \\
99-rc-c  & [\textbf{32}, \textbf{64}]   & 597.67 & [20,  108]   & 2532.76 & [32, 64] \\
101-rc  & [\textbf{36}, 112]   & 532.65 & [26,  144]   & 505.96 & [36, 76] \\
103-rc-c  & [\textbf{32}, \textbf{80}]   & 1419.11 & [22,  120]   & 3400.80 & [32, 80] \\
105-rc-e  & [\textbf{33}, \textbf{112}]   & 1282.90 & [21,  124]   & 410.79 & [33, 112] \\
108-rc-e  & [\textbf{39}, \textbf{120}]   & 957.98 & [24,  136]   & 1566.66 & [39, 120] \\
114-rc-c  & [\textbf{30}, \textbf{64}]   & 733.92 & [22,  96]   & 299.50 & [30, 64] \\
121-rc-e  & [\textbf{34}, 116]   & 1034.38 & [24,  128]   & 3100.10 & [34, 104] \\
126-rc-e  & [\textbf{41}, \textbf{120}]   & 2562.32 & [29,  148]   & 2626.65 & [41, 120] \\
126-r-c  & [\textbf{37}, \textbf{116}]   & 1485.59 & [24,  144]   & 807.31 & [37, 116] \\
144-rc-c  & [\textbf{34}, \textbf{104}]   & 2325.39 & [22,  136]   & 2332.08 & [34, 104] \\
154-c-c  & [\textbf{15}, \textbf{36}]   & 1719.34 & [6,  68]   & 669.42 & [15, 36] \\
165-r-c  & [\textbf{47}, 220]  & 1614.09 & [34, \textbf{212}]  & 3294.70 & [47, 212] \\
167-r-e  & [\textbf{49}, 204] & 1884.33 & [34,  204]  & 1667.48 & [49, 196] \\
173-r-c  & [\textbf{43}, -] & - & [32,  196]  & 2657.03 & [43, 192] \\
173-rc-c  & [\textbf{39}, \textbf{116}]  & 2955.97 & [24,  164]  & 3203.29 & [39, 116] \\
181-r-e  & [\textbf{54}, \textbf{204}]  & 3349.71 & [35,  \textbf{204}]  & 2369.20 & [54, 204] \\
185-c-c  & [\textbf{20}, \textbf{48}]  & 1216.37 & [12,  60]  & 2561.48 & [20, 48] \\
187-rc-e  & [\textbf{48}, \textbf{128}]  & 2645.47 & [32,  192]  & 2310.65 & [48, 128] \\
198-c-c  & [\textbf{16}, \textbf{36}]  & 487.28 & [8,  68]  & 118.12 & [16, 36] \\
200-r-e  & [\textbf{52}, 288]  & 2545.74 & [32,  \textbf{216}]  & 2152.81 & [52, 216] \\
 \hline
\end{tabular}}
\begin{tablenotes}
\item[c] CPU Intel Core i7 12700F - 4x3.6 GHz, 8x4.9 GHz, 20 threads; RAM 32 GB; OR-Tools CP-SAT 9.6; 3600 sec {time limit}
\end{tablenotes}
\end{threeparttable}	
}
\end{center}}
\end{table}

\section{Conclusions}\label{conc}
In this paper, we have discussed several advances for the Parallel Drone Scheduling Traveling Salesman Problem with cooperative drones. In particular, we have proposed a Constraint Programming model coupled with a valid inequality that allows us to find improved lower and upper bounds for the instances proposed in the literature.
Additionally, we demonstrated that the proposed valid inequality can be used to enhance the performance of other methods such as MILP models.

We have also extended the problem into the new Parallel Drone Scheduling Vehicle Routing Problem with cooperative drones, where several trucks are available. For this new extension, we have proposed two alternative Constraint Programming models and a Mixed Integer Programming model. Experimental results suggest that Constraint Programming guarantees better performance, but seems to have scaling issues on large instances, leaving room for future studies on heuristic approaches tailored to the problem.

\section*{Acknowledgements}
The authors are grateful to Minh Ho\`{a}ng Ha and Minh Anh Nguyen for the useful discussions and suggestions, and for having provided the values of the optimized travel times for the drones.

\bibliographystyle{plain}

\end{document}